\documentclass[11pt]{article}  
\usepackage[margin=0.7in]{geometry}
\usepackage{url}
\usepackage[utf8]{inputenc}
\usepackage[english]{babel}
\usepackage{amsthm}

\usepackage[mathscr]{euscript} 
 \let\mathscr\relax 
\usepackage{cite,mathtools} 
\usepackage{latexsym}
\usepackage{amssymb,amsbsy,amsmath,amsfonts,amssymb,amscd,amsfonts,mathrsfs}
\usepackage{graphicx}
\usepackage{color} 
\usepackage{float}
\usepackage{txfonts}
\usepackage{mathtools}

\newtheorem{remark}{Remark}[]

\numberwithin{equation}{section}

\usepackage{titling}
\predate{\begin{flushleft}}
\postdate{\end{flushleft}}
\allowdisplaybreaks
\usepackage[font=rm, labelfont=bf]{caption}

\title{The impact of phenotypic heterogeneity on chemotactic self-organisation}

\author{\normalsize
 Fiona R Macfarlane$^{1,*}$,
 Tommaso Lorenzi$^{2}$, 
 Kevin J Painter$^{3}$
}

\date{ \small \hspace{3em}  $^{1}$ School of Mathematics and Statistics, University of St Andrews, St Andrews, Scotland;\\
 \hspace{3em}  $^{2}$ Department of Mathematical Sciences ``G. L. Lagrange'', Politecnico di Torino, Torino, Italy;\\
  \hspace{3em}  $^{3}$ Inter-university Department of Regional and Urban Studies and Planning, Politecnico di Torino, Torino, Italy;\\
 \hspace{3em} {\footnotesize $^*$Corresponding authors: frm3@st-andrews.ac.uk; }
 }

\begin{document}
\maketitle

\abstract{The capacity to aggregate through chemosensitive movement forms a paradigm of self-organisation, with examples spanning cellular and animal systems. A basic mechanism assumes a phenotypically homogeneous population that secretes its own attractant, with the well known system introduced more than five decades ago by Keller and Segel proving resolutely popular in modelling studies. The typical assumption of population phenotypic  homogeneity, however, often lies at odds with the heterogeneity of natural systems, where populations may comprise distinct phenotypes that vary according to their chemotactic ability, attractant secretion, {\it etc}. To initiate an understanding into how this diversity can impact on autoaggregation, we propose a simple extension to the classical Keller and Segel model, in which the population is divided into two distinct phenotypes: those performing chemotaxis and those producing attractant. Using a combination of linear stability analysis and numerical simulations, we demonstrate that switching between these phenotypic states alters the capacity of a population to self-aggregate. Further, we show that switching based on the local environment (population density or chemoattractant level) leads to diverse patterning and provides a route through which a population can effectively curb the size and density of an aggregate. We discuss the results in the context of real world examples of chemotactic aggregation, as well as theoretical aspects of the model such as global existence and blow-up of solutions.}


\maketitle

\section{Introduction}\label{sec1}

Chemotaxis describes a motility response in which a cell or organism steers itself according to the concentration gradient of a chemical orienteering signal, termed a chemoattractant (repellent) when movement is up (down) the gradient. At the level of a population, chemotaxis can lead to self-organisation, driving a dispersed population into aggregated groups. This process relies on self-reinforcement, for example a population in which each member secretes its own attractant. Such ``autocrine signalling'' allows attraction between near neighbours, leading to localised accumulations that steadily grow in size. Numerous examples of this collective organisation can be cited, spanning the natural world from microscopic to macroscopic: in microbiology, both bacteria and slime molds use chemotaxis to form multicellular aggregation mounds \cite{budrene1991,budrene1995,bonner2009}; during embryonic development, chemotactic self-organisation in the mesenchyme is a key component of hair and feather morphogenesis in the skin \cite{lin2009,glover2017,bailleul2019,ho2019}; during immune responses, activated neutrophils start to produce chemoattractants that recruit other immune cells~\cite{afonso2012ltb4,lammermann2013neutrophil,glaser2021positive}; in the animal world, marine invertebrates such as starfish \cite{hall2017} and sea cucumbers \cite{marquet2018} release water-borne factors to attract neighbours, a prelude to mass spawning. For further examples, we refer the reader to \cite{painter2019}. 

\smallskip
A large modelling literature has, naturally, emerged. At a population level, the majority of models rely on the partial differential equation (PDE) framework introduced by Keller \& Segel \cite{keller1970}, and a ``minimal'' model can be formulated with straightforward functional forms and featuring a single (we assume cellular) population and its chemoattractant, see Figure \ref{fig:schematic}(a). Following a non-dimensionalisation (see Appendix \ref{app_minimal}), this minimal model is given by the following PDE system
\begin{equation}
\begin{cases}
\dfrac{\partial n}{\partial t}=\nabla\cdot \left(D \nabla n-\chi n \nabla s \right),\\\\
\dfrac{\partial s}{\partial t}=\nabla^2 s + n - s,
\end{cases}
\quad {\bf x}\in\Omega.
\label{minimal_PDE}
\end{equation}
Here, the real, non-negative functions $n \equiv n({\bf x},t)$ and $s \equiv s({\bf x},t)$ represent, respectively, the density of cells and the concentration of chemoattractant at time $t \in \mathbb{R}_{+}$ and at position ${\bf x} \in \Omega$. The set $\Omega$ is an open and bounded subset of $\mathbb{R}^d$ with smooth boundary $\partial \Omega$ and $d \geq 1$ depending on the biological problem under study.

\smallskip
The population dynamics are governed by diffusion (modelling random movement) and chemotactic advection according to the gradient of the attractant, while the chemoattractant is produced by the population, decays and diffuses. For an initially dispersed (i.e. spatially uniform) population to aggregate requires that $\chi > D$. This effectively stipulates that the positive feedback process in which a population secretes its own attractant must overcome the homogenising tendencies of diffusion. Models based on \eqref{minimal_PDE} have received significant interest for their patterning and mathematical properties, see the reviews in \cite{horstmann2003,hillen2009,bellomo2015,painter2019}. Mathematical interest has particularly centred on the question of global existence or blow-up (singularity formation), where for the formulation \eqref{minimal_PDE} blow-up has been shown to occur in the biologically relevant case of $d \geq 2$.
While indicative of self-aggregation, the formation of singularities is problematic and artificial in the context of applications, where one naturally expects the growth of an aggregate to be curbed, for example due to lack of space. Consequently, various plausible ``regularisations'' have been proposed (for example see~ \cite{hillen2009}).

\smallskip
An implicit assumption underlying~\eqref{minimal_PDE} and related models is phenotypic homogeneity of the cell population: members are assumed, as a first approximation,
identical with regards to their ability to perform chemotaxis, grow, produce/secrete attractant, {\em etc}. This, however, ignores the typical phenotypic heterogeneity observed within populations: for example, a spectrum of phenotypes whereby population members express different behaviours to different degrees. Chemotactic heterogeneity has been shown for {\em E. coli} populations which, when subjected to an attractant-filled ``T-maze'' structure, spatially sort according to chemotactic ability \cite{salek2019}. Performing chemotaxis, however, demands energy expenditure and the need to balance energy expenditure against energy intake leads to the natural conclusion that high activity in one area must be compensated by lower activity elsewhere, i.e. a trade-off \cite{keegstra2022ecological}. The negative relationship between chemotaxis and growth investment observed for bacteria \cite{ni2020} provide one such trade-off, but given that numerous functions demand significant energy -- movement, growth, protein synthesis, internal molecule and organelle transport, signal transduction, {\em etc}. -- a wide variety of trade-offs will be necessary \cite{keegstra2022ecological}. 

\smallskip
Intra-population phenotypic heterogeneity can also arise through anti-cooperative behaviours, exemplified by  ``cheater'' strains within the slime mold {\em D. discoideum}. These ``social amoebae'' \cite{bonner2009} are famous for the multicellular phase to their lifecycle, starvation inducing auto-aggregation in response to the chemoattractant cAMP, ultimately manifesting in a differentiated fruiting body in which the majority of the population are preserved as ``spores'', but at the cost of a sacrificed 20\% in the supporting ``stalk''. Cheater strains aim to minimise their stalk contribution, feasibly through a spectrum of strategies that include limiting or avoiding the energy-consuming production of cAMP necessary for aggregation \cite{shaulsky2007}. Intra-population phenotypic heterogeneity also extends to large organisms, where the heterogeneity extends from natural variability to distinctly different sub-populations, for example males and females. In the context of chemotactic heterogeneity, for example, only male sea cucumbers release pheromone attractants during the aggregation leading to spawning \cite{marquet2018}.

\smallskip
Motivated by the examples cited above, in this paper we address the following question: How does phenotypic heterogeneity impact on autoaggregation of a population? To this end, we consider a simple extension of the system \eqref{minimal_PDE} whereby the population is subdivided into distinct chemotactic and secreting phenotypes, with the possibility of switching between them. Following the formulation of the model (Section \ref{sec2}), we use linear stability analysis to derive conditions under which patterning is possible (Section \ref{linear_stab}), and carry out numerical simulations to determine the patterning properties (Section \ref{sec3}). We conclude with a discussion (Section \ref{discussion}), in the process proposing some future model extensions and speculating on the analytical properties of the model in the context of global existence and blow-up of solutions.

\section{Model formulation}\label{sec2}

\begin{figure}[t!]
    \centering
    \includegraphics[width=\textwidth]{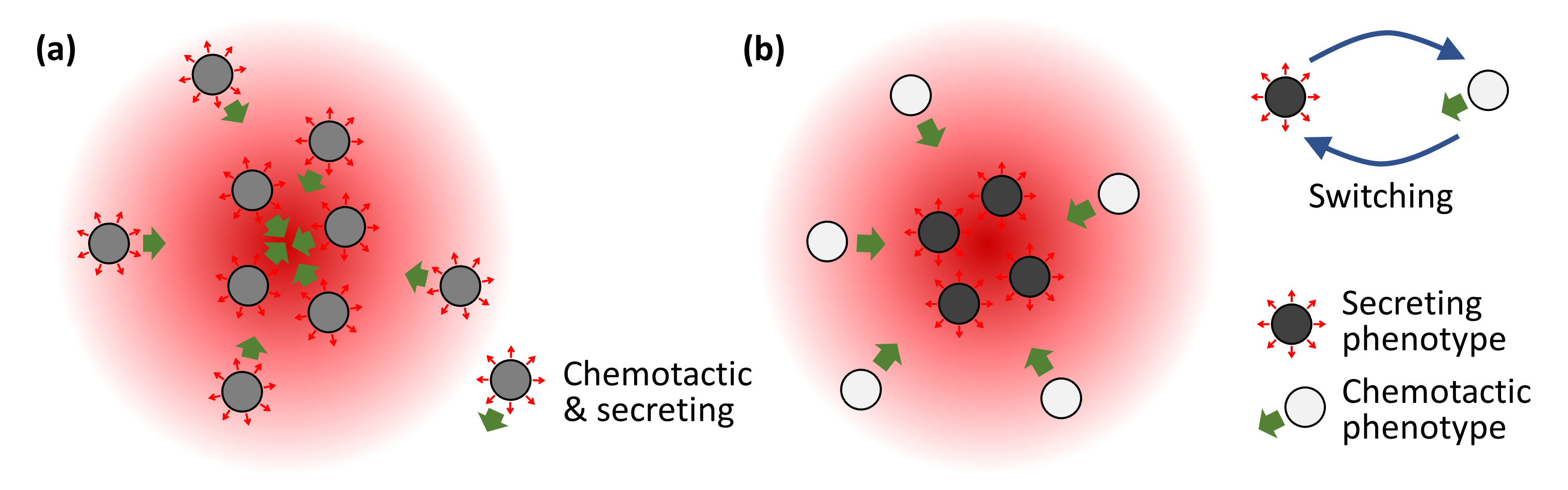}
    \caption{(a) Assumptions underlying the classical minimal autoaggregation model, whereby a population comprises identical members that are both chemosensitive and secrete the chemoattractant signal. (b) Assumptions underlying our binary phenotype model. Here the population is subdivided into distinct chemotactic and secreting phenotypes, with the possibility of switching between these two phenotypic states.}
    \label{fig:schematic}
\end{figure}

We investigate an extension of the basic model~\eqref{minimal_PDE} such that, based on the idea of a trade-off between chemotactic sensitivity and attractant secretion, a binary variable $p \in \{0, 1\}$ is introduced to describe the cell phenotypic state: cells in the phenotypic state $p=0$ produce the chemoattractant but have negligible chemotactic sensitivity, whereas cells in the phenotypic state $p=1$ do not produce the chemoattractant but are chemotactically sensitive, see Figure \ref{fig:schematic}(b). 
As a result, denoting the chemotactic sensitivity of cells in the phenotypic state $p$ as $\chi_{p}$ and the rate of attractant production of cells in the phenotypic state $p$ as $\alpha_{p}$, we assume
$$
\chi_0 = 0, \quad \chi_1 \in \mathbb{R}^*_+, \quad \alpha_0 \in \mathbb{R}^*_+, \quad \alpha_1 = 0\, ,
$$
where $\mathbb{R}^*_+$ is the set of positive real numbers. 

\smallskip
The dynamics of the density of cells in the phenotypic state $p \in \{0, 1\}$ at time $t \geq 0$, $n_p \equiv n_{p}({\bf x},t)$, and the local concentration of chemoattractant, $s \equiv s({\bf x},t)$, are governed by the following system of partial differential equations (PDEs)
\begin{equation}
\begin{cases}
\dfrac{\partial n_0}{\partial t} =D_n\nabla^2 n_0 -{\gamma_{01}(\rho,s)}  \, n_0 + {\gamma_{10}(\rho,s)} \, n_1,\\\\
\dfrac{\partial n_1}{\partial t} =\nabla \cdot \left(D_n \nabla n_1-\chi_1 n_1 \nabla s \right) + {\gamma_{01}(\rho,s)} \, n_0-{\gamma_{10}(\rho,s)} \, n_1, \\\\
{\rho({\bf x},t) := n_0({\bf x},t) + n_1({\bf x},t)}, \\\\
\dfrac{\partial s}{\partial t} =D_s\nabla^2 s + \alpha_0 \, n_0 - \eta \, s,
\end{cases}
\quad {\bf x}\in\Omega.
\label{original_PDE}
\end{equation}
Here, $D_n \in \mathbb{R}^*_+$ is the rate of random motion of the cells, $D_s \in \mathbb{R}^*_+$ is the diffusion rate of the chemoattractant, and $\eta \in \mathbb{R}^*_+$ is the decay rate of the chemoattractant.
{{Generally, eukaryotic cells utilise membrane extensions (pseudopods) to move, whereby the pseudopods anchor the cell to the substrate (or other cells) and permit traction. In the absence of signals promoting directed motion, these pseudopods typically extend in random directions~\cite{van2004chemotaxis} and the inherent random motion is often described using diffusion terms in mathematical models of chemotaxing populations~\cite{painter2019}. To account for this, we have assumed that while cells in the phenotypic state $p=0$ do not undergo directed motion, they still undergo a degree of random motion and are not  completely stationary. For convenience, we assume this inherent degree of randomness is the same for the two phenotypic states, however this could be relaxed at the cost of an additional parameter. }}
 The additional kinetic terms on the right-hand sides of the PDEs~\eqref{original_PDE} for $n_0$ and $n_1$ describe potential switching between phenotypes, that is, the functions
$$
{\gamma_{01}} : \mathbb{R}^2_+ \to \mathbb{R}^*_+ \quad \text{and} \quad {\gamma_{10}} : \mathbb{R}^2_+ \to \mathbb{R}^*_+
$$
represent the rate of switching from phenotypic state 0 (state 1) to phenotypic state 1 (state 0). A switching dependence according to the total cell density $\rho$ could result from a cell changing phenotype following direct contacts at the cell surface, or mediated through a secreted quorum-sensing factor \cite{striednig2022} that provides cell density information. To maximise the generality of the modelling framework, we do not make specific assumptions on these functions and will utilise classical forms based on Hill functions, thus ensuring that switching rates remain bounded.

\smallskip
For the analysis and numerical simulations that follow, we will assume zero-flux conditions at the boundary $\partial \Omega$, i.e. 
\begin{equation}
{\bf{u}}\cdot \nabla n_0=0, \quad {\bf{u}} \cdot \nabla n_1 = 0, \quad {\bf{u}} \cdot \nabla s=0,\quad {\bf x}\in\partial \Omega
\label{boundary}
\end{equation} 
where ${\bf{u}}$ is the unit normal to $\partial \Omega$ that points outwards from $\Omega$. Moreover, initially we simply assume 
$$
n_0({\bf x},0)\ge 0, \quad n_1({\bf x},0)\ge 0, \quad  s({\bf x},0) \ge 0 \quad \forall \, {\bf x} \in \Omega.
$$

\smallskip
We remark that, when subject to boundary conditions~\eqref{boundary}, the PDE system~\eqref{original_PDE} is such that the total cell number is conserved, i.e.
\begin{equation}\label{consrho}
\int_{\Omega} \rho({\bf x},t)\,{\rm d}{\bf x} = \int_{\Omega} \rho({\bf x},0)\,{\rm d}{\bf x} \quad \forall \, t \geq 0.
\end{equation} 
Thus, also the mean of the total cell density over $\Omega$ is conserved. This leads to a further notable quantity $\sigma \in \mathbb{R}_+$ with
\begin{equation}\label{defsigma}
\dfrac{1}{\vert \Omega \vert} \int_{\Omega} \rho({\bf x},0)\,{\rm d}{\bf x} =: \sigma, \end{equation} 
where $\vert \Omega \vert$ is the measure of the set $\Omega$.

\subsection{Dimensionless model}
As detailed in Appendix~\ref{app_nondim}, the PDE system~\eqref{original_PDE} subject to boundary conditions~\eqref{boundary} can be rewritten in the following dimensionless form
\begin{equation}
\begin{cases}
\dfrac{\partial n_0}{\partial t} =D\nabla^2 n_0 {{-\mu_{01}(\rho,s) \, n_0 +\mu_{10}(\rho,s) \, n_1,}}\\\\
\dfrac{\partial n_1}{\partial t} =D \nabla^2 n_1 -\chi \nabla\cdot \left(n_1 \nabla s\right){{+\mu_{01}(\rho,s) \, n_0 -\mu_{10}(\rho,s) \, n_1}}, \\\\
{\rho({\bf x},t) := n_0({\bf x},t) + n_1({\bf x},t)}, \\\\
\dfrac{\partial s}{\partial t} =\nabla^2 s+ n_0 -s,
\end{cases}
\quad {\bf x}\in\Omega
\label{nondim_PDE}
\end{equation}
subject to boundary conditions~\eqref{boundary}.

Here, $D \in \mathbb{R}^*_+$ is the dimensionless rate of random motion of the cells, $\chi \in \mathbb{R}^*_+$ is the dimensionless chemotactic sensitivity, and
\begin{equation}
\mu_{01} : \mathbb{R}^2_+ \to \mathbb{R}^*_+ \quad \text{and} \quad \mu_{10} : \mathbb{R}^2_+ \to \mathbb{R}^*_+
\label{assmu}
\end{equation}
are the dimensionless phenotypic switching functions. Furthermore, upon nondimensionalisation, the conservation relation~\eqref{defsigma} reduces to
\begin{equation}\label{defsigmared}
\dfrac{1}{\vert \Omega \vert} \int_{\Omega} \rho({\bf x},0)\,{\rm d}{\bf x} = 1.
\end{equation}

\subsection{Definitions of phenotypic switching functions}
\label{def_switch}
The definitions of the the phenotypic switching functions, $\mu_{01}$ and $\mu_{10}$, considered here are listed in Table~\ref{tab:1}. Case A corresponds to environment-independent phenotypic switching, with an identical (constant) switching rate, $\mu \in \mathbb{R}^*_+$, in either direction (i.e. from phenotypic state 0 to phenotypic state 1, and vice versa). Cases B$_1$ and B$_2$ each describe density-dependent phenotypic switching: in Case B$_1$ the rate of switching from phenotypic state 0 (state 1) to phenotypic state 1 (state 0) increases (decreases) with the total cell density; in Case B$_2$ these relationships are reversed. Cases C$_1$ and C$_2$ are equivalent to B$_1$ and B$_2$, but where the switching rate now depends on the attractant concentration. In Cases B and C, the parameter $\mu \in \mathbb{R}^*_+$ denotes the maximum rate of switching, while the parameter $q \in \mathbb{R}^*_+$ determines the steepness of the step for the Hill function forms that are used to define $\mu_{01}$ and $\mu_{10}$.

\section{Linear stability analysis of positive uniform steady states}
\label{linear_stab}
In this section, we summarise the results of linear stability analysis of the positive uniform steady states of the PDE system~\eqref{nondim_PDE}, complemented with definition
\begin{equation}
{G(n_0,n_1,\rho,s)} := -\mu_{01}(\rho,s) \, n_0 +\mu_{10}(\rho,s) \, n_1, \label{nondim_defineG}
\end{equation}
and subject to boundary conditions~\eqref{boundary} and condition~\eqref{defsigmared}, which we carried out to identify conditions under which patterning (i.e. formation of cell aggregates) is possible.  {{The parameter regimes that lead to aggregations of the cell population in the minimal model~\eqref{minimal_PDE} have been identified through performing linear stability analysis at the homogeneous steady state of the system and this method can also be applied to investigate patterning in similar chemotactic systems~\cite{hillen2009user}. We explore the conditions under which diffusion/chemotaxis-driven patterning is possible, that is, that the positive uniform steady states of~\eqref{nondim_PDE} are stable to small spatially homogeneous perturbations but unstable to small spatially nonhomogeneous perturbations due to the presence of diffusion/chemotaxis~\cite{murray2001mathematical}.}}To distil the essence of the problem, we focus on the case where the spatial domain is a 1D interval of length $L \in \mathbb{R}^*_+$, i.e. $\Omega := (0, L)$ and ${\bf x} \equiv x$. Extension to higher spatial dimensions is straightforward, but adds little in terms of insight into the mechanisms that underpin the formation of cellular aggregates. The full details of the analysis are provided in Appendix~\ref{app_linear}.

\smallskip
When phenotypic switching does not occur (i.e. when $\mu_{01}\equiv 0$ and $\mu_{10}\equiv 0$), positive uniform steady states are of the form
 \begin{equation}\label{steadystate_app_noswitchpap}
 (n^{\star}_0,n^{\star}_1,s^{\star}) = (\overline{n}, 1- \overline{n}, \overline{n}), \quad 0 < \overline{n} = \dfrac{1}{\vert \Omega \vert} \int_{\Omega} n_0(x,0) \, {\rm d}x<1,
\end{equation}
whereas when phenotypic switching occurs (i.e. when assumptions~\eqref{assmu} hold) for all the definitions of the functions $\mu_{01}$ and $\mu_{10}$ given in Table~\ref{tab:1}, the unique positive uniform steady state is
\begin{equation}\label{steadystate_app_switchpap}
(n^{\star}_0,n^{\star}_1,s^{\star}) = (\overline{n}, 1- \overline{n}, \overline{n}), \quad \overline{n} = 0.5.
\end{equation}
In the former case, the positive uniform steady states~\eqref{steadystate_app_noswitchpap} cannot be driven unstable by small spatially nonhomogeneous perturbations and, therefore, cell aggregates are not expected to form. 

\smallskip
On the other hand, in the latter case, we introduce the definitions
\begin{equation}
H_0:=\left[\frac{\partial G}{\partial n_0}\right]_{ss}, \qquad H_1:=\left[\frac{\partial G}{\partial n_1}\right]_{ss}, \qquad H_s:=\left[\frac{\partial G}{\partial s}\right]_{ss},\label{defineH_text}
\end{equation}
where $[]_{ss}$ indicates that the functions inside the square brackets are evaluated at the steady state~\eqref{steadystate_app_switchpap}. Note that, under assumptions~\eqref{assmu}, we have ${H_1 - H_0 > 0}$. 
{{Furthermore, for the positive uniform steady state~\eqref{steadystate_app_switchpap} to be stable with respect to small spatially homogeneous perturbations, we require that
\begin{equation}
H_1-H_0-H_s \ge 0\label{defineH_conditions1}.
\end{equation}}}

Then the conditions for the positive uniform steady state~\eqref{steadystate_app_switchpap} to be driven unstable by small spatially nonhomogeneous perturbations (i.e. for the formation of cell aggregates to occur) depend on the form of the function $H_1$ and the value of the parameter $\chi$. In more detail:
\begin{itemize}
\item if $\mu_{01}$ and $\mu_{10}$ are such that $H_1<0$ then a {{sufficient}} condition for the formation of cell aggregates to occur is that 
\begin{equation}
\chi>\frac{D(H_{1}-H_{0}-H_{s})}{H_{1} (1-\overline{n})},\label{globalcondition}
\end{equation}
then the perturbation modes labelled by the indices $m \in \mathbb{N}$ such that the following condition on the domain size is satisfied 
\begin{equation}
L>\sqrt{\frac{2D m^2 \pi^2}{-\left[(H_{1}-H_{0}) +D\right]+\sqrt{\left[(H_{1}-H_{0}) -D\right]^2+4 H_{1}\chi (1-\overline{n})+4H_{s}D}}}\label{mindomain_General}
\end{equation}
can grow into aggregation patterns.
\item on the other hand, if $\mu_{01}$ and $\mu_{10}$ are such that $H_1>0$, the following condition holds
\begin{equation}
\chi>\frac{\left[ (H_1-H_0-H_s)+(H_1-H_0)^2\right](D+1)+(3D+1)(H_1-H_0)+2D}{(-H_1)(1-\overline{n})}.\label{globalcondition_H1_less}
\end{equation}
{{Note, in this case the conditions on domain size are less easily calculated, and have therefore been omitted.}}
\end{itemize}

Explicit forms of conditions~\eqref{globalcondition}-\eqref{globalcondition_H1_less} for the particular choices of the phenotypic switching functions considered here are summarised in Table~\ref{tab:1}.

\begin{table}[h]
\resizebox{\textwidth}{!}{
\begin{tabular}{ |c|l|l|l| } 
 \hline
 Case & Switching functions &  Conditions on chemotactic sensitivity   & Conditions on domain size  \\ 
 \hline
 & & &\\
 A & $\begin{aligned} \mu_{01}&	\equiv\mu \\ \mu_{10}&	\equiv \mu \end{aligned}$  &  $\displaystyle{\chi > \frac{2D}{(1-\overline{n})}}$ & $\displaystyle{L>\sqrt{\frac{2D m^2 \pi^2}{-(2\mu+D)+\sqrt{(2\mu-D)^2+4\mu \chi (1-\overline{n}) }}}}$\\
 & & &\\
 & & &\\
 B$_1$ & $\begin{aligned} \mu_{01}&:=\frac{\mu \rho^q}{1 +\rho^q} \\ \mu_{10}&:=\frac{\mu }{1 + \rho^q}\\
 \rho&:=n_0+n_1\end{aligned}$  &  $\begin{aligned} &\displaystyle{\chi > \frac{4D}{ (2-q)(1-\overline{n})}} &(q<2)\\ &\displaystyle{\chi > \frac{4\left[\mu^2 (D +1)+4\mu D+2\mu+2D\right]}{\mu(q-2)(1-\overline{n})} } &(q>2)
 \end{aligned}$ & $\begin{aligned} &\displaystyle{L>\sqrt{\frac{2D m^2 \pi^2}{-\left(\mu +D\right)+\sqrt{\left(\mu -D\right)^2+\mu(2-q)\chi (1-\overline{n})}}}}  &(q<2)\\ &\text{Not calculated symbolically{*}}   &(q>2)
 \end{aligned} $\\
 & &  &\\
   & &  &\\
 B$_2$ & $\begin{aligned} \mu_{01}&:=\frac{\mu }{1 +\rho^q} \\ \mu_{10}&:=\frac{\mu \rho^q}{1 + \rho^q}\\
 \rho&:=n_0+n_1\end{aligned}$ &  $\displaystyle{\chi > \frac{4 D}{(q+2) (1-\overline{n})}}$ & $\displaystyle{L>\sqrt{\frac{2D \pi^2}{-\left(\mu +D\right)+\sqrt{\left(\mu -D\right)^2+ \mu(q+2)\chi (1-\overline{n})}}}}$\\
 & &  &\\
       &  & &\\
 C$_1$& $\begin{aligned} \mu_{01}&:=\frac{\mu s^q}{\overline{n}^q + s^q} \\ \mu_{10}&:=\frac{\mu \overline{n}^q}{\overline{n}^q + s^q} \end{aligned}$  &  $\displaystyle{\chi > \frac{D\left(4\overline{n}+q\right)}{ 2\overline{n} (1-\overline{n})}}$ & $\displaystyle{L>\sqrt{\frac{2D m^2 \pi^2}{-\left(\mu +D\right)+\sqrt{\left(\mu -D\right)^2+2\mu \chi (1-\overline{n})-\frac{\mu q}{\overline{n}}D}}}}$\\
& &  &\\
  & &  &\\
  C$_2$ & $\begin{aligned} \mu_{01}&:=\frac{\mu \overline{n}^q}{\overline{n}^q + s^q} \\ \mu_{10}&:=\frac{\mu s^q}{\overline{n}^q + s^q} \end{aligned}$  &  $\displaystyle{\chi > \frac{D\left(4\overline{n}-q\right)}{ 2\overline{n} (1-\overline{n})}}  \hspace{5em} (q\le4\overline{n}) $ & $\displaystyle{L>\sqrt{\frac{2D m^2 \pi^2}{-\left(\mu +D\right)+\sqrt{\left(\mu -D\right)^2+2\mu\chi (1-\overline{n})+\frac{\mu q}{\overline{n}}D}}}}\hspace{5em}(q\le4\overline{n})$ \\
 & &  &\\
 \hline
\end{tabular}}
\caption{Different forms of phenotypic switching functions investigated along with the resulting conditions~\eqref{globalcondition}-\eqref{mindomain_General}. Note that here $\mu \in \mathbb{R}^*_+$ (i.e. assumptions~\eqref{assmu} hold) and, according to~\eqref{steadystate_app_switchpap}, $\overline{n}=0.5$. {Note that for the specific parameter setting used to produce the results shown in Figure~\ref{1D_compare_oscillations} for the phenotypic switching case $B_1$ we have that as $q$ increases, the minimum domain length that permits patterns $L$ also increases.}}
\label{tab:1}
\end{table}

\paragraph{Remark}
Commenting on the constant switching function case (Case A in Table~\ref{tab:1}) for its simplicity, we note that the condition placed on the chemotactic sensitivity is more stringent than under the minimal formulation of the model (see Section~\ref{sec1}). This is logical, given the division into separate chemotacting and secreting phenotypes: population members actively performing chemotaxis do not reinforce the signal until transitioning into the secreting phenotype, at which point they are no longer actively climbing the attractant gradient. Notably, the rate of phenotypic transition does not impact on the minimal chemotactic sensitivity condition, but it can alter whether self-aggregation will occur through the domain size condition: slow transitioning phenotypes will demand larger domains for pattern formation to occur. The above observations generally hold for the more complicated switching functions, although greater subtleties can arise.

\section{Numerical simulations}\label{sec3}

We carry out numerical simulations of the PDE system~\eqref{nondim_PDE}, complemented with definition~\eqref{nondim_defineG} and subject to boundary conditions~\eqref{boundary} and condition~\eqref{defsigmared}. We consider both one- and two-dimensional spatial domains. The full details of the initial conditions and the set-up of numerical simulations are provided in Appendix~\ref{app:numerics_setup}. We investigate the various forms for the phenotypic switching functions, $\mu_{01}$ and $\mu_{10}$, listed in Table~\ref{tab:1}.  

\subsection{Autoaggregation in 1D}
\begin{figure}[t!]
\centering
\includegraphics[width=\textwidth]{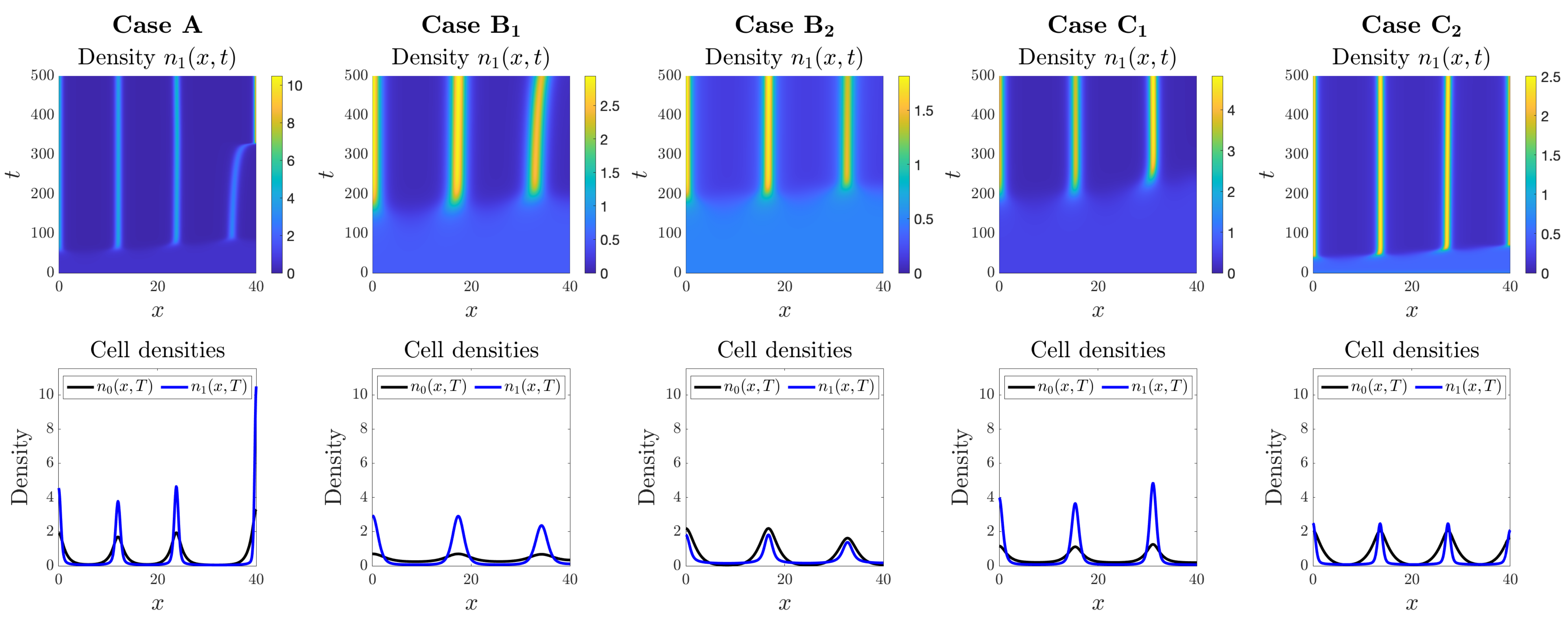}
\caption{{\bf{Comparing the effects of different phenotypic switching functions in 1D.}} {\bf Top row:} Each panel displays the spatial distributions of the cell density $n_1$ over time. {\bf Bottom row:} Each panel displays the spatial distributions of the cell densities $n_0$ (black) and $n_1$ (blue) at the end of simulations, i.e. at $t=T$ with $T:=500$. Each column displays the results of simulations carried out using the definitions of the functions $\mu_{01}$ and $\mu_{10}$ given in Table~\ref{tab:1}, with $\mu=1$ and $q=1$. The values of the parameter $\chi$ used to perform simulations are such that the conditions for pattern formation provided in Table~\ref{tab:1} are satisfied, that is,  $\chi=$ 10 (Case A), 15 (Case B$_1$), 5 (Case B$_2$), 10 (Case C$_1$), 10 (Case C$_2$). The full details of the initial conditions and numerical simulation set-up are provided in Appendix~\ref{app:numerics_setup_1D}.}
\label{1D_compare_switching}
\end{figure}

We first consider a one-dimensional setting, i.e. ${\bf x} \equiv x$, and study the patterns that emerge under each choice of the functions $\mu_{01}$ and $\mu_{10}$ defined in Table~\ref{tab:1}. The corresponding results of numerical simulations are summarised in Figure~\ref{1D_compare_switching}, where the top row of panels displays the plot of the cell density $n_1(x,t)$, for the same choice of initial conditions and parameter values, except for the functions $\mu_{01}$ and $\mu_{10}$ and the chemotactic sensitivity, $\chi$. In the bottom row of panels we display the corresponding cell densities after a simulation time $t=500$, i.e. $n_0(x,500)$ and $n_1(x,500)$. The value of $\chi$ is chosen in each case so as to ensure that the conditions for pattern formation listed in Table~\ref{tab:1} are satisfied. 

\smallskip
In all cases we observe the formation of cell density aggregations, where the height and width of peaks varies with the choice of phenotypic switching functions. These results support the idea that phenotypic switching naturally curbs the density of aggregates. Generally, the chemotactic phenotype forms a concentrated group at the core of the aggregate, while the secreting phenotype is more diffusely spread about the centre.

\begin{figure}[t!]
\centering
\includegraphics[width=\textwidth]{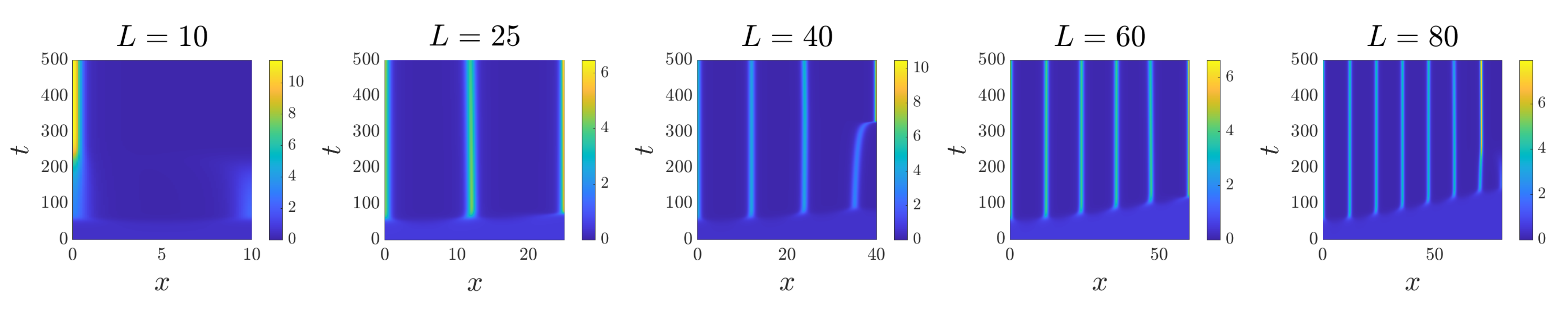}\\
\vspace{-0.5em}
\includegraphics[width=\textwidth]{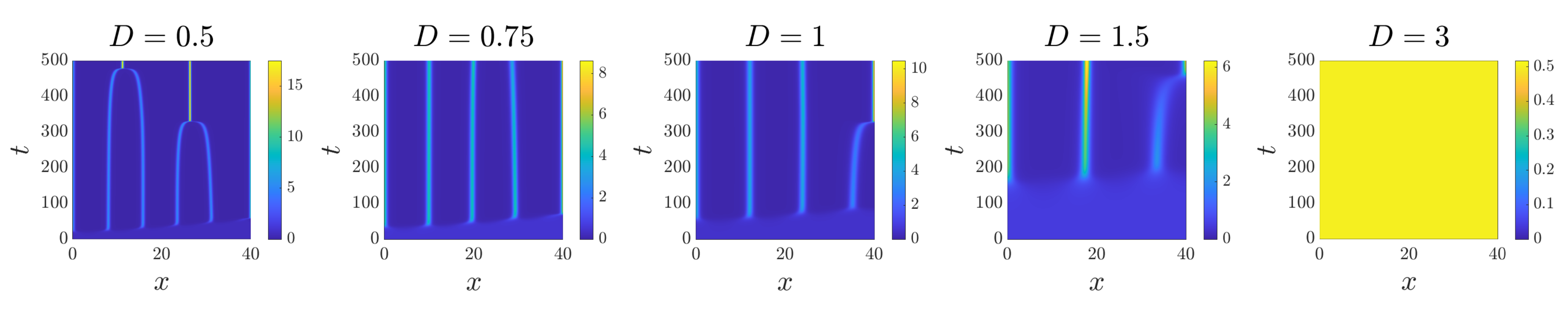}\\
\vspace{-0.5em}
\includegraphics[width=\textwidth]{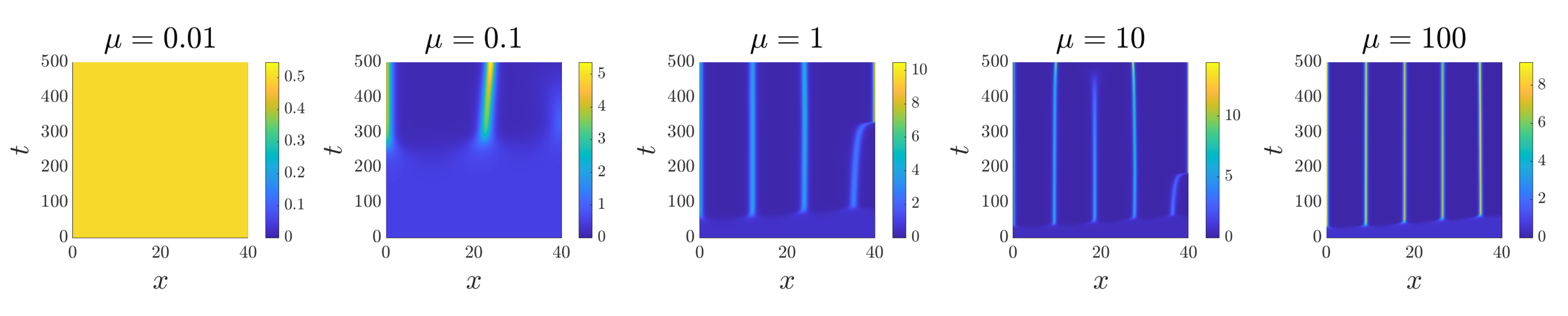}\\
\vspace{-0.5em}
\includegraphics[width=\textwidth]{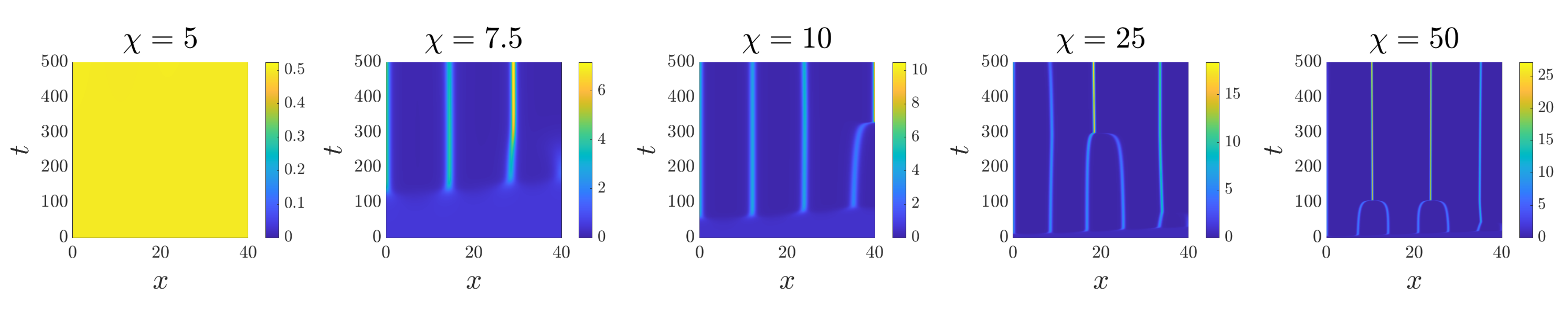}
\caption{{\bf{ Investigating the role of the model parameters in 1D.}} Each panel displays the spatial distribution of the cell density $n_1$ over time. The functions $\mu_{01}$ and $\mu_{10}$ are taken to be those of Case A in Table~\ref{tab:1}, that is, $\mu_{01}\equiv\mu_{10}\equiv\mu$. For each panel we have the general parameter setting $L=40$, $D=1$, $\mu=1$, and $\chi=10$, while the value of one of these parameters is varied in each case as highlighted by the panel titles. The full details of the initial conditions and numerical simulation set-up are provided in Appendix~\ref{app:numerics_setup_1D}.}
\label{1D_constant_test_parameters}
\end{figure}

\smallskip
To further investigate the role of model parameters, we consider the constant phenotypic switching case (i.e. Case A in Table~\ref{tab:1}) and perform simulations while varying one parameter at a time; the results are displayed in Figure~\ref{1D_constant_test_parameters}. This confirms the results of linear stability analysis summarised in Section~\ref{linear_stab}, in particular the instability condition and minimum domain length provided in Table~\ref{tab:1} for Case A. As expected from the characteristic wavelength associated with the fastest growing mode, increasing the domain length $L$ leads to the emergence of a larger number of aggregates, Figure~\ref{1D_constant_test_parameters} (top row). Increasing the random cell movement parameter $D$ both reduces the number of peaks that form and increases the time taken for patterns to emerge from the uniform initial distribution of cells, Figure~\ref{1D_constant_test_parameters} (second row). Too large a value of  $D$ leads to the elimination of pattern formation, since the value of $D$ determines the counter-aggregation diffusion and thereby influences the minimal condition for $\chi$ to have pattern formation. Although the instability condition on $\chi$ is independent of the value of $\mu$, altering $\mu$ does impact on the minimum domain length and, consequently, we observe a loss of pattern formation when the rate of switching becomes too low, Figure~\ref{1D_constant_test_parameters} (third row). Changing the rate of phenotypic switching also impacts on the timescale of pattern formation, with slower switching leading to slower evolution towards the aggregated state. Finally, we investigate the role of $\chi$, where from the linear stability analysis we expect pattern formation beyond a critical $\chi$. This is again confirmed by the simulation results, Figure~\ref{1D_constant_test_parameters} (bottom row). We note that while the linear stability analysis predicts the dynamics when close to the uniform steady state, it becomes less relevant when patterns have formed. Here, numerical simulations indicate dynamics in which nearby aggregates{{, consisting of both cell phenotypes,}} merge leading over time to a reduction in the number of aggregates. This behaviour is well known within chemotaxis-type models, for example see \cite{Potapov2005}. Overall, the results of Figure~\ref{1D_constant_test_parameters} corroborate the conditions for patterning established via linear stability analysis.

\subsection{Oscillating patterns and extinction scenarios in 1D}

\begin{figure}[t!]
\centering
\includegraphics[width=\textwidth]{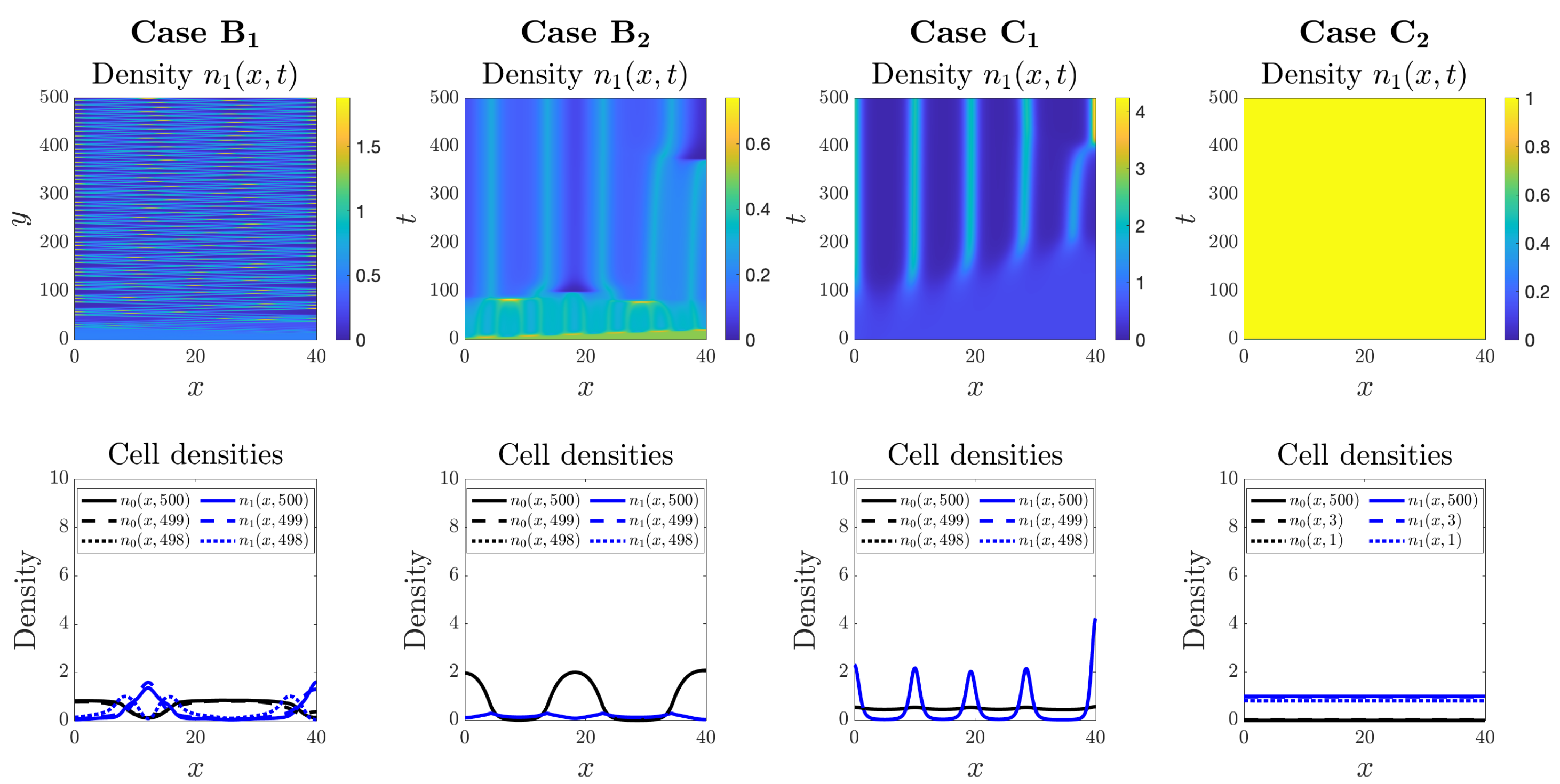}
\caption{{\bf{Comparing the effects of different step-like phenotypic switching functions in 1D.}} {\bf Top row:} Each panel displays the spatial distribution of the cell density $n_1$ over time. {\bf Bottom row:} Each panel displays the spatial distributions of the cell densities $n_0$ (black) and $n_1$ (blue) at three simulation times, i.e. $t=498$, $t=499$, and $t=500$ {{or $t=1$, $t=3$, and $t=500$}}. Each column displays the results of simulations carried out using the definitions of the functions $\mu_{01}$ and $\mu_{10}$ given in Table~\ref{tab:1} (Cases B and C), with $\mu=1$ and $q=30$. The values of the parameter $\chi$ used to perform simulations are such that the conditions for pattern formation provided in Table~\ref{tab:1} are satisfied, that is, $\chi$ is equal to 15 (Case B$_1$), 5 (Case B$_2$), 75 (Case C$_1$), 10 (Case C$_2$). The full details of the initial conditions and numerical simulation set-up are provided in Appendix~\ref{app:numerics_setup_1D}. The full time evolution of the cell densities for each case are provided in Supplementary Material video SM1.}
\label{1D_compare_oscillations}
\end{figure}

We investigate further the cases of environment-dependent phenotypic switching functions given in Table~\ref{tab:1} for step-like changes from a negligible rate of switching to a maximum rate of switching, as some threshold in total cell density (Cases B) or attractant concentration (Cases C) is breached. Specifically, we consider the same general parameter setting as for Figure~\ref{1D_compare_switching}, but with $q$ increased (from $q=1$ to $q=30$) to generate the step-like change and $\chi$ chosen to satisfy conditions for instability provided in Table~\ref{tab:1}. Results are displayed in Figure~\ref{1D_compare_oscillations}.

\smallskip
Introducing step-like switching with respect to the total population density (Cases B) can profoundly impact on the form and nature of aggregates. For the setting in which high total cell densities trigger a switch from the chemotactic to secreting state (Case B$_2$), we observe a significant flattening of the peaks: the maximum density within each aggregate decreases and the peak broadens. This suggests that controlling the balance of secreting and chemotactic cells according to the total population level can be used as an effective means of balancing the density of aggregates. In the reverse setting, where high total densities trigger a switch from the secreting state to the chemotactic state (Case B$_1$), we find evidence of novel patterning. Specifically, non-stationary patterning in which {{complementary aggregates of each individual cell phenotype}} undergo (apparently) sustained oscillations in time. Spatio-temporal oscillations often arise in situations in which the corresponding linear stability analysis indicates the presence of eigenvalues with both a positive real part and non-zero imaginary components (for example see~\cite{tania2012role}). To understand whether this is the case here, we calculate the eigenvalues for this setting. Specifically, we calculate the eigenvalues, i.e. the roots of~\eqref{characteristic_EQ}, for various values of $\mu$ and $\chi$, under Case B$_1$ with other parameters as listed in Figure~\ref{1D_compare_oscillations}. As shown in Figure~\ref{eigenvalues_oscillations}, as the values of $\chi$ and $\mu$ are altered we do indeed observe complex eigenvalues that switch from a negative to positive real part.

\smallskip
Introducing a step-like switching according to the chemoattractant concentration (Cases C) reveals the possibility of further dynamics. For the setting in which higher attractant concentrations trigger a switch from secreting to chemotactic (Case C$_1$), self-organisation occurs as previously, but again with relatively low aggregate densities and a secreting phenotype spread almost uniformly through space. Case C$_2$, where a higher attractant concentration triggers a switch from chemotactic to secreting, demonstrates the possibility of extinction scenarios: {{rapid}} evolution to a population dominated by the chemotaxis phenotype. Here a feedback loop is formed, in which as the number of cells with secreting phenotype drops, so does the overall level of attractant and phenotypic switching is increasingly
weighted towards the chemotactic state.

\subsection{Emergence of patterns in 2D}
{{To consider a more biologically relevant situation, we now turn to a two dimensional setting, i.e. ${\bf x} \equiv (x,y)$.  We note that the transition into two spatial dimensions raises immediate (analytical) questions regarding global existence and blow-up of solutions: for the corresponding minimal model \eqref{minimal_PDEnd}, blow-up typically occurs in two dimensions for scenarios in which autoaggregation is predicted (e.g. see the reviews \cite{horstmann2003,hillen2009,bellomo2015}). In the 2D setting, we again investigate the emergence of patterns under each definition of the phenotypic switching functions listed in Table~\ref{tab:1}.} Results are displayed in Figure~\ref{2D_compare_switching}, where in each column from top to bottom we plot the spatial distributions of $n_0$ and $n_1$ for each phenotypic switching case at the same time-step ($t=500$). Parameters are as in Figure~\ref{1D_compare_switching}.

\smallskip
As predicted from the linear stability analysis, if parameters are selected to satisfy the instability conditions, we observe pattern formation. Specifically, we see the emergence of a spot-like pattern of aggregates, although the size and sharpness of aggregates varies with the form of phenotypic switching. Comparing the distribution of the densities of cells in the two phenotypic states, we observe that the distribution of the chemotactic state is generally concentrated in a sharper peak at the core of each aggregate, with the secreting state more dispersed about the centre.

\smallskip
In comparison to the constant switching form, total cell density dependent and chemoattractant dependent switching forms lead to reduced densities within cell aggregates, clearest within Case B$_1$. Notably, none of the switching cases that have been considered here lead to numerical blow-up phenomena, defined as instances in which the numerical solutions form densities and/or gradients that lead to numerical instability and simulation failure. This hints that the introduction of switching between chemotactic and secreting phenotypic states may lead to global existence of solutions, although caution is noted given the numerical nature of the study. We return to this in the discussion, where we also exploit radial symmetry scenarios to perform more refined simulations.

\begin{figure}[t!]
\centering
\includegraphics[width=\textwidth]{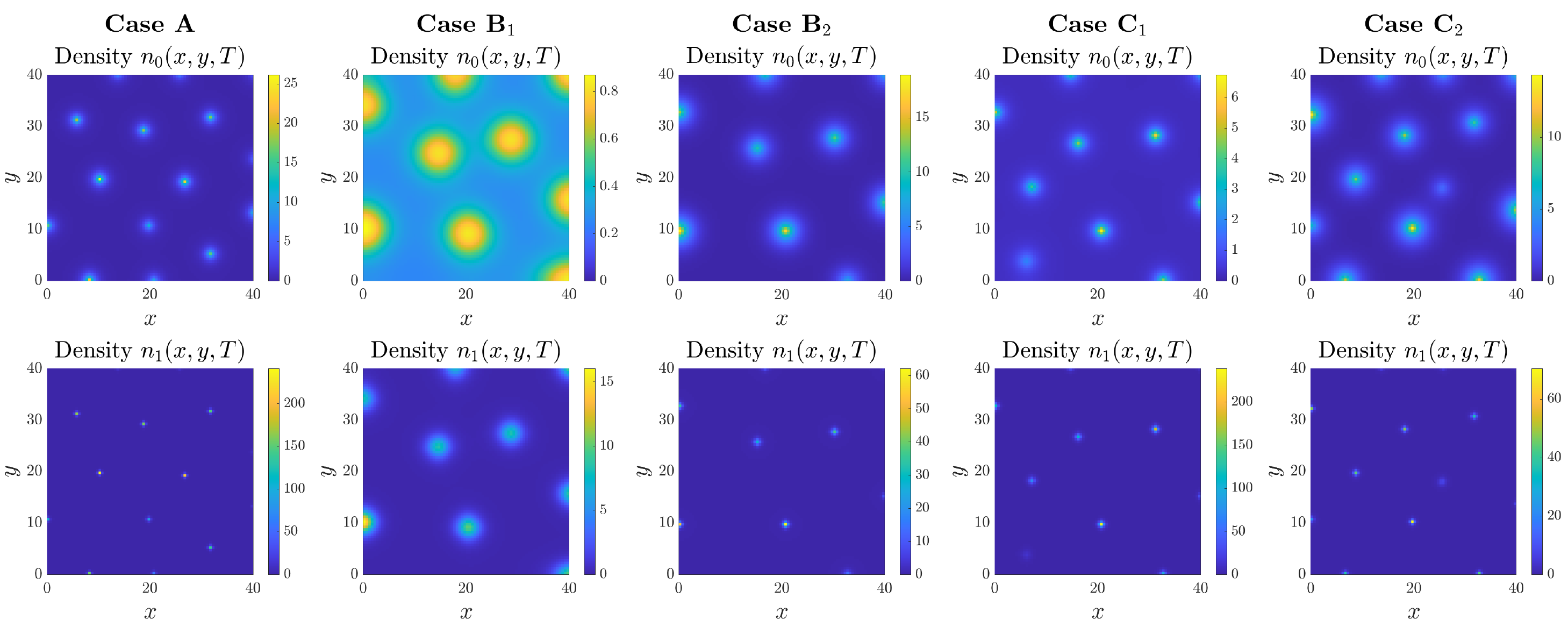}
\caption{{\bf{ Comparing the effects of different phenotypic switching functions in 2D.}} {\bf Top row:} Each panel displays the 2D spatial distribution of the cell density $n_0$ at the time $t=T$ with $T:=500$. {\bf Bottom row:} Each panel displays the 2D spatial distribution of the cell density $n_1$ at the time $t=T$ with $T:=500$. Each column displays the results of simulations carried out using the definitions of the functions $\mu_{01}$ and $\mu_{10}$ given in Table~\ref{tab:1}, with $\mu=1$ and $q=1$. The values of the parameter $\chi$ used to perform simulations are such that the conditions for pattern formation provided in Table~\ref{tab:1} are satisfied, that is, $\chi$ is equal to 10 (Case A), 15 (Case B$_1$), 5 (Case B$_2$), 10 (Case C$_1$), 10 (Case C$_2$). The full details of the initial conditions and numerical simulation set-up are provided in Appendix~\ref{app:numerics_setup_2D}. {{The full time evolution corresponding to each panel is provided in Supplementary Material video SM2.}}}
\label{2D_compare_switching}
\end{figure}

\smallskip
We also consider two dimensional equivalents for the results displayed in  Figure~\ref{1D_compare_oscillations}, i.e. where step-like phenotypic switching functions were selected. The results are shown in {{Figure~\ref{2D_compare_oscillations}}}, where again each column displays the density of the two phenotypes, $n_1$ and $n_2$, at $t=500$. {{For comparison, we additionally plot the density of $n_1$ at $t=490$, i.e. the density at an earlier time instant.} Similar phenomena to those observed in one dimension are observed. For Case B$_1$ we see a two-dimensional analogue to the oscillating pattern, with aggregated {{single cell phenotype}} structures that undergo sustained temporal dynamics. For Case B$_2$ we observe the evolution to relatively low-density cell aggregates, in which the peak density at the core of the aggregate remains bounded at a relatively low level with respect to the uniform density. In this scenario there is a notable distribution of the two phenotypes, with the secreting phenotype centered at the core of an aggregate and encapsulated by a chemotaxing ring at the periphery. This is distinct from the standard arrangement (chemotactic phenotype concentrated in a high density peak at the centre) and arises through a transition from chemotaxis to secretion behaviour as an individual reaches the higher density at the centre of an aggregate. For Case C$_1$ we see the emergence of low-density aggregates where, as in the one-dimensional setting, the secreting phenotype is spread almost uniformly in space while the chemotactic phenotype forms spot-like patterns.
\begin{figure}[t!]
\centering
\includegraphics[width=0.8\textwidth]{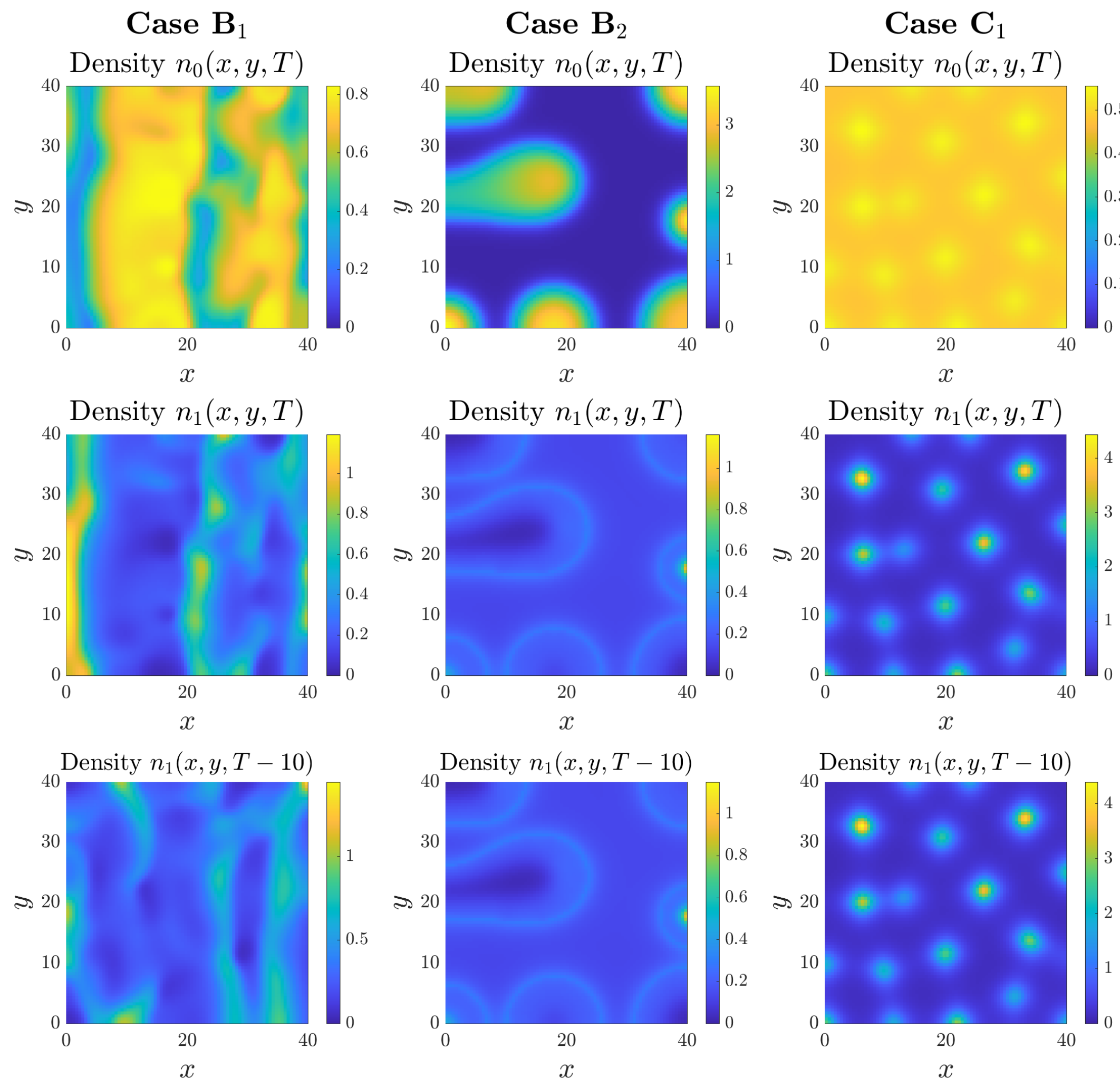}
\caption{{\bf{Comparing the effects of different step-like phenotypic switching functions in 2D.}} {\bf Top row:} Each panel displays the 2D spatial distribution of the cell density $n_0$ at the time $t=T$ with $T:=500$. {\bf Middle row:} Each panel displays the 2D spatial distribution of the cell density $n_1$ at the time $t=T$. {\bf Bottom row:} Each panel displays the 2D spatial distribution of the cell density $n_1$ at a previous time instant, i.e. $t=T-10=490$. Each column displays the results of simulations carried out using the definitions of the functions $\mu_{01}$ and $\mu_{10}$ given in Table~\ref{tab:1} (Cases B and C), with $\mu=1$ and {{$q=30$}}. The values of the parameter $\chi$ used to perform simulations are such that the conditions for pattern formation provided in Table~\ref{tab:1} are satisfied, that is, $\chi$ equal to 15 (Case B$_1$), 5 (Case B$_2$), 75 (Case C$_1$). The full details of the initial conditions and numerical simulation set-up are provided in Appendix~\ref{app:numerics_setup_2D}. {{The full time evolution corresponding to each panel is provided in Supplementary Material video SM3.}}}
\label{2D_compare_oscillations}
\end{figure}

\smallskip
The two sets of simulation results under Case B$_2$ (displayed in Figure~\ref{2D_compare_switching} and Figure~\ref{2D_compare_oscillations}) suggest that this form of phenotypic switching is particularly effective for controlling the density and width of cell aggregates: for example, leading to densities only a few times larger than the uniform steady state density, as compared to several hundred times under constant switching forms (Figure~\ref{2D_compare_switching}). We investigate this further, considering a range of parameter settings for Case B$_2$. Results are shown in Figure~\ref{2D_compare_densitydep}, where each panel displays the {{density $n_1$}} obtained under a distinct parameter setting (cf. panel titles). {{Note that we have omitted here the figures for the corresponding densities $n_0$, where in all cases the aggregates consist of cells of both phenotypes.}} The top row explores the impact of increasing the chemotactic sensitivity, where we observe a transition from low-density stripe-like patterns at low $\chi$, to ring-like structures with a peaked centre for moderate values of $\chi$, and to sharp high-density peaks at higher $\chi$. The transitions under increasing $\mu$ (middle row) or $q$ (bottom row) show this patterning in reverse, with the transition from high-density spots to lower-density rings and/or stripes as the parameters are increased. Overall, these results indicate that density-dependent phenotypic switching functions can lead to a broad spectrum of spatial patterning, with regimes of spatio-temporal patterning for certain functional forms, or allow for a spectrum between stripe-like and spot-like patterns.

\begin{figure}[t!]
\centering
\includegraphics[width=0.8\textwidth]{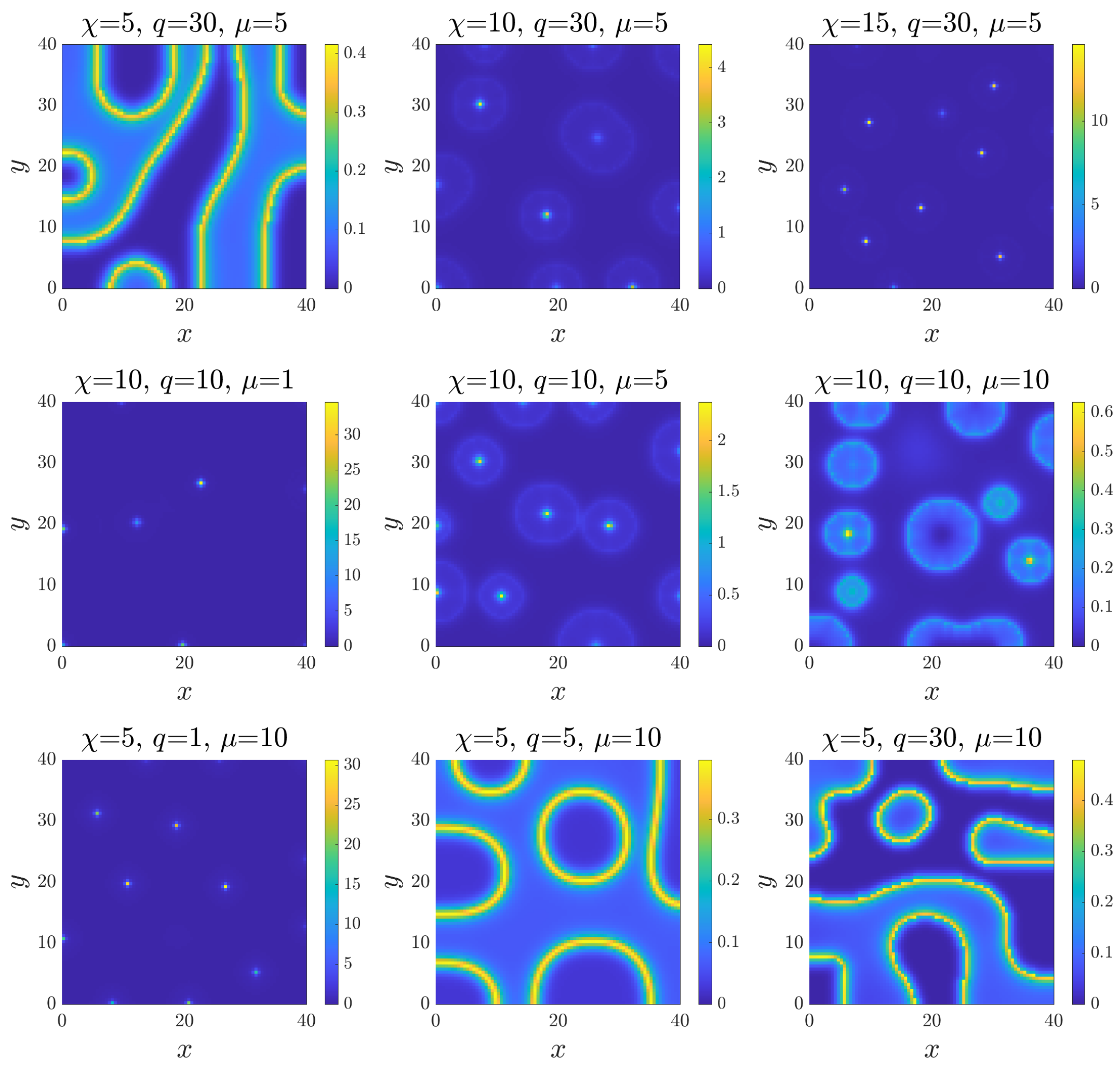}
\caption{{\bf{Investigating the role of the model parameters in 2D.}} Each panel displays the 2D spatial distribution of the cell density $n_1$ at the time $t=T$ with $T:=500$. The functions $\mu_{01}$ and $\mu_{10}$ are taken to be those of Case B$_2$ in Table~\ref{tab:1}. For each panel, the parameter setting investigated is highlighted by the panel titles. The full details of the initial conditions and numerical simulation set-up are provided in Appendix~\ref{app:numerics_setup_2D}. {{The full time evolution corresponding to each panel is provided in Supplementary Material video SM4.}}}
\label{2D_compare_densitydep}
\end{figure}

\section{Discussion and conclusions} \label{discussion}
\paragraph{Overview of results}
We have developed a simple model to describe the dynamics of a phenotypically heterogeneous population of cells that interacts with some chemoattractant. Intra-population phenotypic heterogeneity is incorporated into the model by assuming the population to be composed of individuals in two distinct phenotypic states, with the possibility of switching between these states. Individuals in one phenotypic state perform chemotaxis but do not secrete attractant, whilst individuals in the other phenotypic state secrete attractant but do not display chemotaxis. 

\smallskip
Using a combination of linear stability analysis and numerical simulations, we have shown that pattern formation can occur in the form of cell density aggregates that contain a mixture of the two phenotypic states. Moreover, when the rate of switching between these phenotypic states depends on the local environment -- the total cell density or the attractant concentration -- a range of dynamical behaviours can be observed, including oscillations in space and time, pattern structures that range from stripes to spots, or ``phenotype extinction'' scenarios in which the cell population evolves to one in which all members occupy the same phenotypic state.

\paragraph{Discussion of results in the context of chemotactic aggregation}
Notably, our results indicate that a level of transitioning between phenotypic states is an essential criterion for autoaggregation within the system: when switching is negligible, aggregation is not possible. Certain systems, however, may preclude this transitioning: for example, an ecological population featuring a division into male and female members. In the context of aggregation scenarios, many marine invertebrates employ chemical aggregation to achieve higher population densities prior to broadcast swarming (simultaneous release of gametes into the water column). Recent studies in certain sea cucumbers suggest that only the males release the aggregating pheromones, but that these attract nearby members of both sexes \cite{marquet2018}. In this respect, the modelling here suggests that male to male attraction is essential for aggregation to occur in this system.

\smallskip
The inherent cost attached to some labour (movement, protein synthesis, {\em etc.}) leads naturally to the question of energy balance, i.e. that individuals must balance their energy expenditure with the rate at which energy can be generated. Intra-population phenotypic diversity may then results from a trade-off: individuals in a phenotypic state that allows them to perform one task at a high level have less energy for other tasks \cite{keegstra2022ecological}. The trade-off considered here is between performing chemotaxis and both producing and secreting an attractant. In particular, we have shown that a division of labour between the time spent in the phenotypic state enabling individuals to perform chemotaxis and the time spent in the phenotypic state enabling individuals to produce and secrete the attractant can still lead to successful aggregation of the population, provided there is a non-negligible degree of transition between these phenotypic states. 
{{Here, we have concentrated on energy trade-offs related to the chemoattractant secretion and chemotactic movement of cells. However, we should note that the approach taken here has been phenomenological in spirit and a plausible future extension would be to explicitly account for energetics. One way to do this would be to initially formulate a model at the discrete level, for example through a random walk for an individual cell, but in which its movement and secretion behaviour depend on an internal variable that represents its energy reserve. Accounting for a dependency on the state of an internal variable has a rich history in the modelling of chemotaxis behaviour (e.g.~\cite{erban2004}) and, using similar methods that we have employed before~\cite{bubba2020discrete,macfarlane2020hybrid,chaplain2020bridging,macfarlane2022individual}, we could then derive the continuum limit as a PDE model that can further be analysed.
We considered here the trade-off to be between cell motility and chemoattractant secretion, however, other trade-offs could exist in the population such as proliferation-migrations trade-offs. Such trade-offs have been observed in cancer growth, that is the ``go-or-grow'' hypothesis, and have been studied extensively both experimentally and theoretically~\cite{corcoran2003,gallaher2019impact,giese1996dichotomy,hatzikirou2012go,hoek2008vivo,pham2012density,stepien2018,vittadello2020,zhigun2018,stinner2016global,dhruv2013reciprocal,giese2003cost,wang2012ephb2,xie2014targeting}. In fact, mathematical models of the form developed in this manuscript can also be used to study these proliferation-migration trade-offs~\cite{macfarlane2022individual,lorenzi2022trade}. In this work, we considered that cells in the phenotypic state $p=0$ would still undergo a degree of random motion, however, an alternative model would be to consider a phenotypic transition into a sessile population (i.e. with neither random or directed movement terms). An investigation into the effects that this would have on the solutions to the system would be of interest.}}

\smallskip
Acquiring an aggregated state can be essential for many populations, whether as a prelude to tissue morphogenesis during embryonic development or as a key stage within the life-cycle of some unicellular or multicellular organisms. Yet aggregating can also be risky: within an ecological context, as an example, reaching a high density state could lead to the depletion of local resources and starvation. Consequently, while a mechanism for aggregating may be crucial, it may be similarly important to have some counter-aggregation mechanism in place, preventing the population from over-accumulating at a location. An upper bound to the density of individuals can certainly arise from considerations of ``volume-filling'' (i.e. individuals being unable to occupy the same space) \cite{hillen2009}, but for many aggregates the density lies far below such extreme conditions, see, for instance, \cite{mittal2003}. Our results support the idea that allowing individuals in the population to transition between chemotactic and secreting phenotypes according to the local environment provides an effective means of controlling the density of aggregates. The greatest control was found in the case of responses to total density, specifically whereby individuals transition from performing chemotaxis to performing secretion as the local density surpasses some threshold. Under this scenario, low density aggregates are readily achieved, as well as the capacity to arrange into a diversity of two-dimensional structures from spots to labyrinthine stripes. 

\paragraph{Discussion of results in the context of PDE models for chemotaxis}
PDE models for chemotaxis that comprise two or more populations have been extensively explored by numerous authors, from both applied and theoretical perspectives. Avoiding a detailed discussion, we restrict to highlighting a few particularly relevant examples. From an applied perspective, the study presented in~\cite{tania2012role} features two populations of forager and scrounger type, where foragers follow (and consume) a nutrient, and scroungers follow the foragers. ``Behavioural switching'' analogous to the population transition terms considered here were employed to describe shifting between these two behaviour types and, despite the absence of an autocrine relationship between the foragers and the attractant, spatio-temporal oscillations were observed. Two population systems have also been considered in the context of ``autocrine-paracrine'' relationships, where two populations can produce and respond to their own attractant (autocrine), but also respond to the attractant produced by the other population (paracrine): see \cite{painter2009} in the context of differential-chemotaxis induced sorting, or \cite{knutsdottir2014} in the context of tumour-macrophage signalling.

\smallskip
From a more theoretical perspective, the model~\eqref{original_PDE} sits closely to the general ``$p$-populations, $q$-chemicals'' models that were described in \cite{wolansky2002}: the addition in the model here lies in the possibility to switch between populations. A substantial literature has emerged on the analytical properties, particularly regarding blow-up of solutions, for example see \cite{espejo2012,tello2012,kurganov2014,fu2016,wang2017time}. These studies raise the question of blow-up in this system, and we briefly comment on the numerically-obtained insights from the present study. For the minimal model~\eqref{minimal_PDE}, blow-up is dimensionally-dependent: in one space dimension, solutions exist globally, but in two or more dimensions finite-time blow-up can occur (e.g. see the reviews in \cite{horstmann2003,hillen2009,bellomo2015}). The studies conducted here, however, suggest that switching between secreting and chemotactic states curbs the rate at which aggregates accumulate, leading to the question as to whether the model permits global existence. To explore this in further depth, we exploit the  convenience of a radially symmetric (two dimensions) setting, thereby reducing to an effectively one-dimensional structure for numerical simulations, which permits computation for a highly refined spatial discretisation (see  Appendix~\ref{app:numerics_setup_radial} for details). Concentrating on the case of constant rate of phenotypic switching (i.e. Case A given in Table~\ref{tab:1} whereby $\mu_{01}\equiv\mu_{10}\equiv\mu$ with $\mu\in \mathbb{R}^*_+$), we track the cell densities over time, computing until there is negligible change to the solution profile or, if relevant, numerical blow-up occurs. The results are displayed in Figure~\ref{2D_radial}, under a range of values of $\mu$. In all cases solutions are found to evolve to smooth profiles and no numerical blow-up was observed. Increasing $\mu$ generates a sharper profile, however the maximum density remains bounded. A similar boundedness to the computed solutions was found to occur in a spherically symmetric (three dimensional) setting (data not shown). While exercising caution given the numerical nature of this study, these results suggest that the model studied here may admit globally existing non-uniform solutions.

\begin{figure}[H]
\centering
\includegraphics[width=\textwidth]{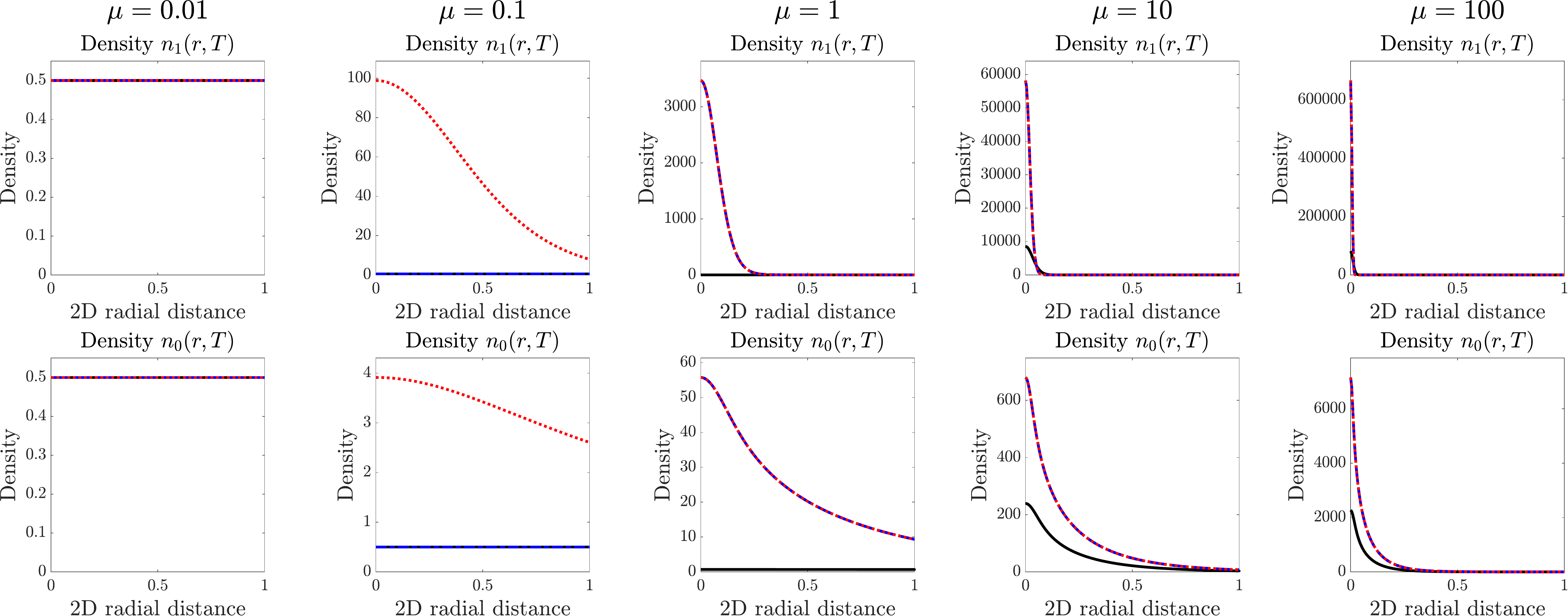}\\
\includegraphics[width=\textwidth]{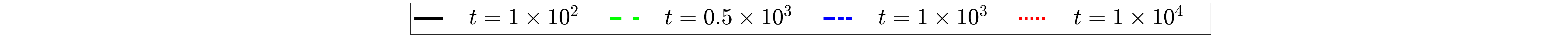}
\caption{{\bf{Investigating the possibility of having no numerical blow-up in the radially symmetric case.}} {\bf Top row:} Each panel displays the 2D radial distance of the cell density $n_1$ at four increasing time instants, as described by the legend. {\bf Bottom row:} Each panel displays the 2D radial distance of the cell density $n_0$ at four increasing time instants, as described by the legend. In each case, phenotypic switching is of the form of Case A given in Table~\ref{tab:1}, and each column displays the simulation results obtained for different value of $\mu \in \mathbb{R}^*_+$. The full details of the initial conditions and numerical simulation set-up are provided in Appendix~\ref{app:numerics_setup_radial}.}
\label{2D_radial}
\end{figure}

\paragraph{Possible generalisations of the model}
The mathematical model defined by the PDE system~\eqref{original_PDE} subject to the boundary conditions~\eqref{boundary} relies on the assumption that the cell phenotypic state is binary. Hence, as a first generalisation of this model, we could consider the case in which the population is subdivided into a discrete spectrum of $P+1$ phenotypes, labelled by the index $p \in \mathcal{P} := \{ 0, \ldots, \cal{P}\}$, which vary according to both their chemotactic ability and their rate of attractant production. The generalised model would then comprise the following PDE system
\begin{equation}
\begin{cases}
\displaystyle{\dfrac{\partial n_p}{\partial t}=\nabla\cdot \left(D_n \nabla n_p - \chi_p n_p \nabla s \right) + \sum_{q\in{\cal{P}}, q\ne p} \Big(\gamma_{qp}(\rho,s) n_q - \gamma_{pq}(\rho,s) n_p \Big)},\\\\
\displaystyle{\rho({\bf x},t) := \sum_{p\in{\cal{P}}} n_p({\bf x},t)}, \\\\
\displaystyle{\dfrac{\partial s}{\partial t} =D_s\nabla^2 s + \sum_{p\in{\cal{P}}} \alpha_p n_p - \eta \, s},
\end{cases}
 {\bf x}\in\Omega,
\label{phenotype_PDE}
\end{equation}
subject to the following zero-flux boundary conditions
$$
{\bf{u}} \cdot \nabla n_p =0, \quad {\bf{u}} \cdot \nabla s=0,\quad x\in\partial \Omega, \quad p=0, \ldots, \cal{P}.
$$
In analogy with the case of model~\eqref{original_PDE}, the function $\gamma_{qp}(\rho, s)$, with $\gamma_{qp} : \mathbb{R}^2_+ \to \mathbb{R}_+$, models the rate at which cells switch from phenotype $p$ to phenotype $q \neq p$.

As a further generalisation, in the vein of~\cite{macfarlane2022individual,lorenzi2022trade,lorenzi2021invasion}, we could also consider the case where the cell phenotypic state varies along a continuum and is thus described by a continuous variable $\theta \in \mathbb{R}$. In this case, the dynamics of the density of cells in the phenotypic state $\theta$ at time $t \geq 0$, $n \equiv n(\theta,{\bf x},t)$, and the local concentration of the chemoattractant would be governed by the following system of partial integro-differential equations
\begin{equation}
\begin{cases}
\displaystyle{\dfrac{\partial n}{\partial t}=\nabla_{{\bf x}}\cdot \left(D_n \nabla_{{\bf x}} n - \chi(y) n \nabla_{{\bf x}} s \right)} \\
\displaystyle{\phantom{.......}+ \int_{\mathbb{R}} \gamma(\theta,\theta',\rho,s) \, n(\theta',{\bf x},t) \, {\rm d}\theta' - n \int_{\mathbb{R}} \gamma(\theta',\theta,\rho,s) \, {\rm d}\theta'},\\\\
\displaystyle{\rho \equiv \rho({\bf x},t) := \int_{\mathbb{R}} n(\theta,{\bf x},t) \, {\rm d}\theta},
\\\\
\displaystyle{\dfrac{\partial s}{\partial t} =D_s\nabla^2_{{\bf x}} s + \int_{\mathbb{R}} \alpha(\theta) \, n(\theta,{\bf x},t) \, {\rm d}\theta - \eta \, s},
\end{cases}
\; {\bf x}\in\Omega,
\label{contphenotype_PDE}
\end{equation}
subject to the following zero-flux boundary conditions
$$
{\bf{u}} \cdot \nabla_{{\bf x}} n = 0, \quad {\bf{u}} \cdot \nabla_{{\bf x}} s=0,\quad {\bf x} \in \partial \Omega.
$$
Here, the functions $\chi : \mathbb{R} \to \mathbb{R}^*_+$ and $\alpha : \mathbb{R} \to \mathbb{R}^*_+$ model, respectively, the chemotactic sensitivity and the rate of attractant production of cells in the phenotypic state $\theta$. Moreover, the function $\gamma(\theta, \theta',\rho,s)$, with $\gamma : \mathbb{R}^2 \times \mathbb{R}^2_+ \to \mathbb{R}_+$, models the rate of switching from the phenotypic state $\theta'$ to the phenotypic state $\theta$.

Analysis and numerical simulation of the generalised models~\eqref{phenotype_PDE} and~\eqref{contphenotype_PDE}, which are expected to pose a series of analytical and numerical challenges that make these models interesting mathematical objects {\em per se}, will allow for further investigation into how phenotypic heterogeneity shapes the emergence of chemotactic self-organisation in different biological contexts.

\paragraph{Acknowledgments}
{\small F.R.M. gratefully acknowledges support from the RSE Saltire Early Career Fellowship `Multiscale mathematical modelling of spatial eco-evolutionary cancer dynamics' (Fellowship No. 1879). T.L. gratefully acknowledges support from the Italian Ministry of University and Research (MUR) through the grant ``Dipartimenti di Eccellenza 2018-2022'' (Project no. E11G18000350001) and the PRIN 2020 project (No. 2020JLWP23) ``Integrated Mathematical Approaches to Socio–Epidemiological Dynamics'' (CUP: E15F21005420006). K.J.P. acknowledges ``MIUR-Dipartimento di Eccellenza'' funding to the Dipartimento Interateneo di Scienze, Progetto e Politiche del Territorio (DIST).
}

\appendix

\section{The minimal model}
\label{app_minimal}

In dimensional form, the minimal model is given by
\begin{equation}
\begin{cases}
\dfrac{\partial n}{\partial t}=\nabla\cdot \left(D_n \nabla n- \chi n \nabla s \right),\\\\
\dfrac{\partial s}{\partial t}=D_s \nabla^2 s + \alpha n - \eta s,
\end{cases}
\quad {\bf x}\in\Omega.
\label{minimal_PDEnd}
\end{equation}
The real, non-negative functions $n \equiv n({\bf x},t)$ and $s \equiv s({\bf x},t)$ represent, respectively, the density of cells and the concentration of chemoattractant at time $t \in \mathbb{R}_{+}$ and at position ${\bf x} \in \Omega$. The set $\Omega$ is an open and bounded subset of $\mathbb{R}^d$ with smooth boundary $\partial \Omega$ and $d \geq 1$ depending on the biological problem under study. The parameters $D_n \in \mathbb{R}^*_{+}$ and $D_s\in \mathbb{R}^*_{+}$ denote the cell and chemoattractant diffusion coefficients, respectively. The parameter $\alpha\in \mathbb{R}^*_{+}$ denotes the chemoattractant production rate per cell, while $\eta\in \mathbb{R}^*_{+}$ is the decay rate. The chemotactic sensitivity coefficient, $\chi$, is positive for an attractant response, i.e. $\chi\in \mathbb{R}^*_{+}$. Under lossless boundary conditions (for example, zero flux) it is straightforward to observe that the positive uniform steady state to \eqref{minimal_PDEnd} is given by $(\overline{n}, \alpha\overline{n}/\eta)$, where $\overline{n} \in \mathbb{R}^*_{+}$ represents the mean initial dispersed density of cells, that is,
$$
\overline{n} := \dfrac{1}{\vert \Omega \vert} \int_{\Omega} n({\bf x},0) \, {\rm d}{\bf x}.
$$
A nondimensionalisation according to 
{\small
\begin{equation}
\hat{{\bf x}}:=\sqrt{\frac{\eta}{D_s}} {\bf x},\quad \hat{t}:= \eta t, \quad \hat{n}:=\frac{n}{\overline{n}}, \quad\hat{s}:=\frac{\eta}{\alpha \overline{n}} s, \quad D:= \frac{D_n}{D_s}, \quad \hat{\chi} = \frac{\chi \alpha \overline{n}}{\eta D_s},
\nonumber
\end{equation}}
leads (after dropping the hats) to the system \eqref{minimal_PDE}, with the positive uniform steady state rescaled to $(1,1)$.

\section{Nondimensionalisation of the PDE system~\eqref{original_PDE}}
\label{app_nondim}
Without loss of generality, focussing on the case of a 1D spatial domain (i.e. ${\bf x} \equiv x$), we introduce the following dimensionless variables
{\small
\begin{equation}
\hat{x}:=\frac{x}{X},\quad \hat{t}:=\frac{t}{T}, \quad \hat{s}:=\frac{s}{S}, \quad \hat{n}_0:=\frac{n_0}{N}, \quad \hat{n}_1:=\frac{n_1}{N}, \quad \hat{\rho}:= \frac{\rho}{N} \nonumber
\end{equation}}
along with the following dimensionless forms of the phenotypic switching functions
{\small
\begin{equation}
 \hat{\gamma}_1\equiv\hat{\gamma}_1(\hat{\rho},\hat{s}):=\frac{\gamma_1(\rho,s)}{\Gamma_1}, \quad  \hat{\gamma}_2\equiv\hat{\gamma}_2(\hat{\rho},\hat{s}):=\frac{\gamma_2(\rho,s)}{\Gamma_2}, \nonumber
\end{equation}}
where $X$, $T$, $S$, $N$, $\Gamma_1$, and $\Gamma_2$ are scale factors. Substituting these into the PDE system~\eqref{original_PDE} posed on a a 1D spatial domain yields
{\small
\begin{align}
\label{nondim_step1}
\frac{N}{T}\frac{\partial \hat{n}_0}{\partial \hat{t}}&= D_n\frac{N}{X^2}\frac{\partial^2 \hat{n}_0}{\partial \hat{x}^2}-N \, \Gamma_1  \hat{\gamma}_1\hat{n}_0+N \, \Gamma_2  \hat{\gamma}_2\hat{n}_1\nonumber\\
\frac{N}{T}\frac{\partial \hat{n}_1}{\partial \hat{t}}&=D_n\frac{N}{X^2}\frac{\partial^2 \hat{n}_1}{\partial \hat{x}^2}-\frac{\chi_1 N S}{X^2} \frac{\partial }{\partial \hat{x}}\left( \hat{n}_1 \frac{\partial \hat{s}}{\partial \hat{x}}\right)+N \, \Gamma_1 \hat{\gamma}_1\hat{n}_0 - N \, \Gamma_2  \hat{\gamma_2}\hat{n}_1\nonumber\\
\frac{S}{T}\frac{\partial \hat{s}}{\partial \hat{t}}&=D_s\frac{S}{X^2}\frac{\partial^2 \hat{s} }{\partial \hat{x}^2}+ \alpha_0 N \hat{n}_0 -\eta S \hat{s}.\nonumber
\end{align}}

Next, isolating the time-derivatives, we find
{\small
\begin{align}
\frac{\partial \hat{n}_0}{\partial \hat{t}}&= D_n\frac{T}{X^2}\frac{\partial^2 \hat{n}_0}{\partial \hat{x}^2}- \Gamma_1 T  \hat{\gamma}_1\hat{n}_0+ \Gamma_2 T  \hat{\gamma}_2\hat{n}_1\nonumber\\
\frac{\partial \hat{n}_1}{\partial \hat{t}}&=D_n\frac{T}{X^2}\frac{\partial^2 \hat{n}_1}{\partial \hat{x}^2}-\frac{\chi_1 S T}{X^2} \frac{\partial }{\partial \hat{x}}\left(\hat{n}_1 \frac{\partial \hat{s}}{\partial \hat{x}}\right)+ \Gamma_1 T  \hat{\gamma}_1\hat{n}_0-\Gamma_2 T  \hat{\gamma}_2\hat{n}_1\nonumber\\
\frac{\partial \hat{s}}{\partial \hat{t}}&=D_s\frac{T}{X^2}\frac{\partial^2 \bar{s} }{\partial \hat{x}^2}+ \frac{\alpha_0 N T}{S} \hat{n}_0 -\eta T  \hat{s}.\nonumber
\end{align}
}

Then, we choose
{\small
$$
T:=\frac{1}{\eta},\quad X:=\sqrt{\frac{D_s}{\eta}}, \quad S:= \frac{\alpha_0 N}{\eta},
$$
}
in order to remove the coefficients in the chemoattractant equation. In so doing, we obtain
{\small
\begin{align}
\frac{\partial \hat{n}_0}{\partial \hat{t}}&= \frac{D_n}{D_s}\frac{\partial^2 \hat{n}_0}{\partial \hat{x}^2}- \frac{\Gamma_1}{\eta}   \hat{\gamma}_1\hat{n}_0+ \frac{\Gamma_2}{\eta}   \hat{\gamma}_2\hat{n}_1\nonumber\\
\frac{\partial \hat{n}_1}{\partial \hat{t}}&=\frac{D_n}{D_s}\frac{\partial^2 \hat{n}_1}{\partial \hat{x}^2}-\frac{\chi_1 \alpha_0 N}{\eta  D_s} \frac{\partial }{\partial \hat{x}}\left(\hat{n}_1 \frac{\partial \hat{s}}{\partial \hat{x}}\right)+ \frac{\Gamma_1}{\eta}   \hat{\gamma}_1\hat{n}_0-\frac{\Gamma_2}{\eta} \hat{\gamma}_2\hat{n}_1\nonumber\\
\frac{\partial \hat{s}}{\partial \hat{t}}&=\frac{\partial^2 \hat{s} }{\partial \hat{x}^2}+ \hat{n}_0 - \hat{s}.\nonumber
\end{align}}

Finally, introducing the dimensionless parameters and functions
{\small
$$
D:=\frac{D_n}{D_s}, \quad \chi := \frac{\chi_1 \alpha_0  N}{\eta  D_s}, \quad \mu_{01}:=\frac{\Gamma_1}{\eta}\hat{\gamma}_1, \quad \mu_{10}:=\frac{\Gamma_2}{ \eta}\hat{\gamma}_2$$
}

and then removing the hats we obtain the dimensionless PDE system ~\eqref{nondim_PDE} with rescaled phenotypic switching terms defined via~\eqref{nondim_defineG} as provided in the main text. Note that in the above we have not yet needed to specify the scaling of the cell densities, i.e. set the choice of $N$. Here we can take advantage of the conservation of the total cell number (cf. conservation relation~\eqref{consrho}), choosing $N=\sigma$ with $\sigma$ defined via relation~\eqref{defsigma}, such that, in the framework of the dimensionless PDE system~\eqref{nondim_PDE}, we have condition~\eqref{defsigmared}.

\section{Positive uniform steady states of the PDE system~\eqref{nondim_PDE} and linear stability analysis}
\label{app_linear}
We perform linear stability analysis of the positive uniform steady states of the PDE system~\eqref{nondim_PDE} complemented with definition~\eqref{nondim_defineG} and subject to boundary conditions~\eqref{boundary} and condition~\eqref{defsigmared}.

\subsection{Positive uniform steady states}
When phenotypic switching does not occur (i.e. if $\mu_{01} \equiv 0$ and $\mu_{10} \equiv 0$), the positive uniform steady states $(n^{\star}_0,n^{\star}_1,s^{\star})$ of the PDE system~\eqref{nondim_PDE} complemented with definition~\eqref{nondim_defineG} and subject to the zero-flux boundary conditions~\eqref{boundary} and condition~\eqref{defsigmared} are of the form
 \begin{equation}\label{steadystate_app_noswitch}
 (n^{\star}_0,n^{\star}_1,s^{\star}) = (\overline{n}, 1- \overline{n}, \overline{n}), \quad 0 < \overline{n} = \dfrac{1}{\vert \Omega \vert} \int_{\Omega} n_0(x,0) \, {\rm d}x<1.
\end{equation}

On the other hand, when there is phenotypic switching (i.e. when assumptions~\eqref{assmu} hold), the positive uniform steady states $(n^{\star}_0,n^{\star}_1,s^{\star})$ of the PDE system~\eqref{nondim_PDE} complemented with definition~\eqref{nondim_defineG} and subject to the zero-flux boundary conditions~\eqref{boundary} and condition~\eqref{defsigmared} are of the form
$$
n^{\star}_0=\overline{n}, \quad n^{\star}_1=1-n^{\star}_0, \quad s^{\star}=n^{\star}_0,
$$
where $\overline{n} \in \mathbb{R}^*_+$ is given by the following algebraic equation
$$
\overline{n} = \frac{\mu_{10}(1,\overline{n})}{\mu_{01}(1,\overline{n})+\mu_{10}(1,\overline{n})}.
$$
For all the definitions of $\mu_{01}$ and $\mu_{10}$ given in Table~\ref{tab:1}, this algebraic equation admits a unique solution, that is, $\overline{n} = 0.5$. Hence, there is a unique positive uniform steady state
\begin{equation}\label{steadystate_app_switch}
(n^{\star}_0,n^{\star}_1,s^{\star}) = (\overline{n}, 1- \overline{n}, \overline{n}), \quad \overline{n} = 0.5.
\end{equation}

\subsection{Linear stability analysis of positive uniform steady states}
In order to study the linear stability of the positive uniform steady states to small perturbations, we focus on the 1D case where the spatial domain is a 1D interval of length $L \in \mathbb{R}^*_+$, i.e. $\Omega := (0,L)$. We make the ansatz
\begin{eqnarray}
n_0(x,t) &=& \overline{n} + \tilde{n}_0 \, \exp{(\lambda t)} \, \varphi_k(x), \nonumber \\ 
n_1(x,t) &=& (1-\overline{n}) + \tilde{n}_1 \, \exp{(\lambda t)} \, \varphi_k(x),\nonumber \\ 
s(x,t) &=& \overline{n} + \tilde{s} \, \exp{(\lambda t)} \, \varphi_k(x), \nonumber
\end{eqnarray}
where $\tilde{n}_0, \tilde{n}_1, \tilde{s} \in \mathbb{R}_*$ with $\vert \tilde{n}_0 \vert \ll 1$, $\vert \tilde{n}_1 \vert \ll 1$ and $\vert \tilde{s} \vert \ll 1$, $\lambda \in \mathbb{C}$ and $\{\varphi_k\}_{k \geq 1}$ are the eigenfunctions of the Laplace operator, acting on functions defined on $(0,L)$ and subject to zero Neumann boundary conditions, indexed by the wavenumber $k$, i.e.
\begin{equation}\label{formk}
k  = \dfrac{m \pi}{L}, \;\; m \in \mathbb{N}.
\end{equation}

Linearising the PDE system~\eqref{nondim_PDE} posed on $(0,L)$ about a positive uniform steady state $(\overline{n}, 1- \overline{n}, \overline{n})$ and using the above ansatz yields the following matrix equation
{\small
\begin{equation}
\begin{bmatrix}
\lambda+k^2D-H_{0} & -H_{1} & -H_{s}\\
H_{0} & \lambda+k^2D+H_{1} & -\chi (1-\overline{n})k^2 +H_s\\
-1 & 0 & \lambda+k^2+1
\end{bmatrix}
\begin{bmatrix}
\tilde{n}_0 \\
\tilde{n}_1 \\
\tilde{s}
\end{bmatrix}=0.
\end{equation}
}
For the above matrix equation to admit a solution $(\tilde{n}_0, \tilde{n}_1, \tilde{s}) \in \mathbb{R}^3_*$, we need    
{\small
\begin{equation}
\label{characteristic_EQ}
\lambda^3+A(k^2)\lambda^2+B(k^2)\lambda+C(k^2)=0,
\end{equation}
}
where
{\small
\begin{align}
\label{Characteristic_defineA}
A(k^2)&:=\left(2D+1\right)k^2+\left(H_{1}-H_{0}+1\right),\\
\label{Characteristic_defineB}B(k^2)&:=D \left(D+2\right) k^4 + \left[(H_{1}-H_{0}+2) D  +(H_1-H_0)\right]k^2 + \nonumber \\ & \quad +\left(H_{1}-H_{0}-H_s\right),\\
\label{Characteristic_defineC}C(k^2)& :=D^2 k^6+\left[D(H_{1}-H_{0}) +D^2\right]k^4 + \nonumber \\ & \quad + \left[D(H_{1}-H_{0}-H_s)- H_{1}\chi (1-\overline{n})\right]k^2. 
\end{align}
}
Here,
{\small
\begin{equation}
\label{defineH}
 H_{0}:=\left[\frac{\partial G}{\partial n_0}\right]_{ss} \quad H_{1}:=\left[\frac{\partial G}{\partial n_1}\right]_{ss} \quad \text{and} \quad H_{s}:=\left[\frac{\partial G}{\partial s}\right]_{ss},
\end{equation}
}
where $G \equiv G(n_0,n_1,\rho,s)$ is defined via~\eqref{nondim_defineG}, that is,
$$
G(n_0,n_1,\rho,s) := -\mu_{01}(\rho,s) \, n_0 +\mu_{10}(\rho,s) \, n_1,
$$
and $[]_{ss}$ indicates that the functions inside the square brackets are evaluated at the positive uniform steady state.

\subsubsection{{Stability w.r.t spatially homogeneous perturbations}}
{We see that when small spatially homogeneous perturbations are considered (i.e. when $k^2=0$) we have that the characteristic equation becomes
$$
\lambda^3+A(0)\lambda^2+B(0)\lambda+C(0)=0,
$$
where,
$$
A(0)=\left(H_{1}-H_{0}+1\right), \quad B(0)=\left(H_{1}-H_{0}-H_s\right) \quad \text{and} \quad C(0)=0.
$$
For a cubic polynomial of the form $\lambda^3 + a \lambda^2 + b \lambda + c$, the Routh-Hurwitz criterion ensures that ${\rm Re}(\lambda)<0$ if and only if $a>0$, $b>0$, $c>0$, and $ab-c>0$.
We begin by noting that, since $\rho := n_0 + n_1$,
$$
\frac{\partial G}{\partial n_0}=-\frac{\partial \mu_{01}}{\partial \rho}n_0-\mu_{01}+\frac{\partial \mu_{10}}{\partial \rho}n_1, \quad \frac{\partial G}{\partial n_1}=-\frac{\partial \mu_{01}}{\partial \rho}n_0+\mu_{10}+\frac{\partial \mu_{10}}{\partial \rho}n_1.
$$
Hence,
\begin{equation}
\label{h1mh0pos}
H_1-H_0=\left[\frac{\partial G}{\partial n_1}\right]_{ss}-\left[\frac{\partial G}{\partial n_0}\right]_{ss}=[\mu_{01}+\mu_{10}]_{ss} > 0.
\end{equation}
Therefore $A(0)>0$ and for a positive uniform steady state $(\overline{n}, 1- \overline{n}, \overline{n})$ to be stable with respect to small spatially homogeneous perturbations we require $B(0)>0$, that is,
\begin{equation}
\label{conditionHs}
\left(H_{1}-H_{0}-H_s\right)\ge0.
\end{equation}
}

\subsubsection{Stability w.r.t spatially nonhomogeneous perturbations}
A positive uniform state $(\overline{n}, 1- \overline{n}, \overline{n})$ will be driven unstable by small spatially nonhomogeneous perturbations (i.e. spatial patterns will emerge) if there exists at least one $k^2 \in \mathbb{R}^*_+$ for which ${\rm Re}(\lambda)>0$. For a cubic polynomial of the form $\lambda^3 + a \lambda^2 + b \lambda + c$, the Routh-Hurwitz criterion ensures that ${\rm Re}(\lambda)<0$ if and only if $a>0$, $b>0$, $c>0$, and $ab-c>0$. Therefore, for patterns to emerge we need at least one of the following conditions 
{\small
\begin{equation}
\label{RH_conditions}
A(k^2)>0, \quad B(k^2)>0, \quad C(k^2)>0, \quad A(k^2) B(k^2) - C(k^2)>0
\end{equation}}
not to be satisfied for some $k^2 \in \mathbb{R}^*_+$.

\begin{remark}
From definitions~\eqref{Characteristic_defineA} and \eqref{Characteristic_defineB}, we have that
{\footnotesize
\begin{align}
A(k^2)B(k^2)&=\left[D\left(D+2\right)\left(2D+1\right)\right] k^6\nonumber\\
& \quad + \left[\left(2D+1\right)\left((H_{1}-H_{0}+2) D  +(H_1-H_0)\right)+\left(H_{1}-H_{0}+1\right)\left(D^2+2D\right)\right] k^4 \nonumber\\
& \quad+\left[ \left(H_{1}-H_{0}-H_s\right)\left(2D+1\right)+ \left(H_{1}-H_{0}+1\right)\left((H_{1}-H_{0}+2) D  +(H_1-H_0)\right)\right]k^2\nonumber\\
& \quad+\left[H_{1}-H_{0}+1\right]\left[H_{1}-H_{0}-H_s\right].\nonumber
\end{align}
}
Using this along with definition~\eqref{Characteristic_defineC}, we calculate $A(k^2)B(k^2)-C(k^2)$. A little algebra yields 
{\small
\begin{align}
A(k^2)B(k^2)-C(k^2)&=\left[2D\left(D+1\right)^2\right] k^6 \nonumber\\
&+ \left[(3D+1)(H_1-H_0)(D+1)+2D(2D+1)\right] k^4 \nonumber\\
&+\left\{\left[ (H_1-H_0-H_s)+(H_1-H_0)^2\right](D+1)\right\}k^2\nonumber\\
&+\left[(3D+1)(H_1-H_0)+2D+H_{1}\chi (1-\overline{n})\right]k^2\nonumber\\
&+\left(H_{1}-H_{0}+1\right)\left(H_{1}-H_{0}-H_s\right).\label{defineABmC}
\end{align}
}
\end{remark}

\subsubsection{Case without phenotypic switching}
When phenotypic switching does not occur (i.e. if $\mu_{01} \equiv 0$ and $\mu_{10} \equiv 0$) then 
$$
{G(n_0,n_1,\rho,s)} \equiv 0 \quad \Longrightarrow \quad H_0=H_1=H_s=0.
$$
Under definitions~\eqref{Characteristic_defineA}-\eqref{Characteristic_defineC}, this implies that
$$
A(k^2) := \left(2D+1\right)k^2+1, \ \  B(k^2) :=D \left(D+2\right) k^4 + 2D k^2, \ \  C(k^2) :=D^2 k^6+D^2k^4,
$$
and expression~\eqref{defineABmC} implies that
$$
A(k^2)B(k^2)-C(k^2) =2Dk^2\left[(D+1)^2k^4 +2(D+1)k^2+1\right].
$$
Hence, conditions~\eqref{RH_conditions} will be satisfied for all $k^2 \in \mathbb{R}^*_+$ and, therefore, we do not expect spatial patterns to emerge.

\subsubsection{Case with phenotypic switching}
We now turn to the case where phenotypic switching occurs (i.e. when assumptions~\eqref{assmu} hold).

Under assumptions~\eqref{assmu} and using \eqref{h1mh0pos}, definition~\eqref{Characteristic_defineA} gives
$$
A(k^2) = \left(2D+1\right)k^2+[\mu_{01}+\mu_{10}]_{ss} + 1 > 0 \quad \forall k^2 \in \mathbb{R}^*_+,
$$
that is, $A(k^2)$ will satisfy condition~\eqref{RH_conditions} for all $k^2 \in \mathbb{R}^*_+$. As a result, for pattern formation to occur we need one of the remaining conditions~\eqref{RH_conditions} to be violated. Since calculations are rather tedious, here we will be limiting ourselves to deriving some conditions on the functions $\mu_{01}$ and $\mu_{10}$ and the parameter $\chi$ that are required to have $A(k^2)B(k^2)-C(k^2)<0$ for some $k^2 \in \mathbb{R}^*_+$ or required to have $C(k^2)<0$ for a certain range of values of $k^2$. {{Recall that for the positive uniform steady states to be stable with respect to homogeneous perturbations we require that condition \eqref{conditionHs} holds. Hence, we require the functions $\mu_{01}$ and $\mu_{10}$ to be such that condition \eqref{conditionHs} is satisfied.}} Since $H_1-H_0>0$ (cf. relation~\eqref{h1mh0pos}), under assumption~\eqref{conditionHs} we have that:
\begin{itemize}
\item due to expression~\eqref{defineABmC}, to have $A(k^2)B(k^2)-C(k^2)<0$ for some $k^2 \in \mathbb{R}^*_+$, it is necessary that the functions $\mu_{01}$ and $\mu_{10}$ are such that $H_1<0$ and the parameter $\chi$ must satisfy the following condition
{\small
\begin{equation}
{\small
\chi>\frac{\left[ (H_1-H_0-H_s)+(H_1-H_0)^2\right](D+1)+(3D+1)(H_1-H_0)+2D}{(-H_1)(1-\overline{n})};\label{globalcondition_H1_lessapp}
}
\end{equation}
}
\item due to definition~\eqref{Characteristic_defineC}, if the functions $\mu_{01}$ and $\mu_{10}$ are such that $H_1>0$ and the parameter $\chi$ satisfies the following condition
{\small
\begin{equation}\label{app_globalcondition}
{\small
\chi>\frac{D(H_{1}-H_{0}-H_{s})}{H_{1} (1-\overline{n})},
}
\end{equation}
}

then the condition $C(k^2)<0$ is satisfied for $k^2 \in \mathbb{R}^*_+$ within the following range of unstable modes
{\small
\begin{equation}
\label{conditionsk_General}
0<k^2<\frac{-\left[(H_{1}-H_{0}) +D\right]+\sqrt{\left[(H_{1}-H_{0}) -D\right]^2+4 H_{1}\chi (1-\overline{n})+4H_{s}D}}{2D}.
\end{equation}
}
\end{itemize}
Moreover, substituting~\eqref{formk} into~\eqref{conditionsk_General} we find that, when assumptions~\eqref{conditionHs}, $H_1>0$ and \eqref{app_globalcondition} hold, the unstable modes are those that are labelled by the wave numbers $m \in \mathbb{N}$ such that the following condition on the domain size is satisfied   
{\small
\begin{equation}
\label{app_mindomain_General}
L>\sqrt{\frac{2D m^2 \pi^2}{-\left[(H_{1}-H_{0}) +D\right]+\sqrt{\left[(H_{1}-H_{0}) -D\right]^2+4 H_{1}\chi (1-\overline{n})+4H_{s}D}}}.
\end{equation}
}

The values of $H_0$, $H_1$ and $H_s$, along with the values of the notable quantities $H_1 - H_0$ and $H_1-H_0-H_s$, for the phenotypic switching functions defined in Table~\ref{tab:1} are provided in Table~C1. 

\begin{table}[h!]
\label{tab:2}
\centering
\begin{tabular}{|c|c|c|c|c|c|} 
 \hline
 Case & $H_0$ & $H_1$ & $H_s$ & $H_1-H_0$ & $H_1-H_0-H_s$\\
 \hline
&&&&&\\
 A & $-\mu$ & $\mu$ & $0$ & $2\mu$ & $2\mu$\\
&&&&&\\
 \hline
&&&&&\\
 B$_1$ & $\frac{-\mu(2+q)}{4}$ & $\frac{\mu(2-q)}{4}$ & $0$ & $\mu$ & $\mu$\\
&&&&&\\
 \hline
&&&&&\\
 B$_2$ & $\frac{\mu(q-2)}{4}$ & $\frac{\mu(2+q)}{4}$ & $0$ & $\mu$ & $\mu$\\
&&&&&\\
 \hline
&&&&&\\
 C$_1$ & $-\frac{\mu}{2}$ & $\frac{\mu}{2}$ & $-\frac{\mu q}{4\overline{n}}$ & $\mu$ & $\mu+\frac{\mu q}{4\overline{n}}$\\
&&&&&\\
 \hline
&&&&&\\
 C$_2$ & $-\frac{\mu}{2}$ & $\frac{\mu}{2}$ & $\frac{\mu q}{4\overline{n}}$ & $\mu$ & $\mu-\frac{\mu q}{4\overline{n}}$\\
&&&&&\\
 \hline
 \end{tabular}
 \vspace{0.5em}
 \caption{Values of $H_0$, $H_1$ and $H_s$ defined via~\eqref{defineH}, along with the values of the notable quantities $H_1 - H_0$ and $H_1-H_0-H_s$, for the phenotypic switching functions given in Table~\ref{tab:1}. Here, $\overline{n}=0.5$.} 
 \end{table}

\section{Numerical simulation set-up}
\label{app:numerics_setup}
In this section, we describe the methods and the parameter settings that we used to carry out numerical simulations of the model. All numerical simulations are performed in {\sc{Matlab}}.

\subsection{Set-up of 1D simulations}
\label{app:numerics_setup_1D}
To obtain the results displayed in Figures~\ref{1D_compare_switching}-\ref{1D_compare_oscillations}, we solved numerically the PDE system~\eqref{nondim_PDE} with ${\bf x} \equiv x \in (0,L)$, complemented with definition~\eqref{nondim_defineG}, and subject to boundary conditions~\eqref{boundary} and condition~\eqref{defsigmared}. Unless stated otherwise in the article, we set $L=40$. Numerical simulations were carried out using the {\sc{Matlab}} function {\bf{pdepe}}. A uniform discretisation of the spatial domain with step $\Delta_x=0.1$ was used. Note that {\bf{pdepe}} uses an adaptive time-step to obtain a solution.
As an initial condition, we chose the following perturbed version of the uniform steady state~\eqref{steadystate_app_switch}
\begin{equation}
{\small
n_0(x,0)=\overline{n}, \quad n_1(x,0)=1-\overline{n}, \quad s(x,0) =\overline{n}+0.01\ \exp\left[-A \left(x-\dfrac{L}{2}\right)^2 \right],
\label{initialcondition_1D}
}
\end{equation}
where $\overline{n} = 0.5$ and $A=1\times10^{4}$. We note that, under the initial data defined via~\eqref{initialcondition_1D}, condition~\eqref{defsigmared} is satisfied. Unless stated otherwise in the article, we set $D=1$. All the other parameter values are given in the captions of Figures~\ref{1D_compare_switching}-\ref{1D_compare_oscillations}. 

\subsection{Set-up of 2D simulations}
\label{app:numerics_setup_2D}
To obtain the results displayed in Figures~\ref{2D_compare_switching}-\ref{2D_compare_densitydep}, we solved numerically the PDE system~\eqref{nondim_PDE} with ${\bf x} \equiv (x,y) \in (0,L) \times (0,L)$, complemented with definition~\eqref{nondim_defineG}, and subject to boundary conditions~\eqref{boundary} and condition~\eqref{defsigmared}. Unless stated otherwise in the article, we set $L=40$. Numerical simulations were carried out using the finite volume scheme described in Appendix~\ref{app:2D_scheme}. A uniform discretisation of the spatial domain with steps $\Delta_x=\Delta_y=0.5$ was used along with a uniform discretisation of the time domain with step $\tau=1\times10^{-3}$.

As an initial condition, we chose the following perturbed version of the uniform steady state~\eqref{steadystate_app_switch}
\begin{equation}
{\small
n_0(x,y,0)=\overline{n}, \quad n_1(x,y,0)=1-\overline{n}, \quad s(x,y,0) = \overline{n} + 0.01 R(x,y),
\label{initialcondition_2D}
}
\end{equation}
where $\overline{n} = 0.5$ and the values attained by the function $R$ are random numbers in $[0,1]$. We note that, under the initial data defined via~\eqref{initialcondition_2D}, condition~\eqref{defsigmared} is satisfied. Unless stated otherwise in the article, we set $D=1$. All the other parameter values are given in the captions of Figures~\ref{2D_compare_switching}-\ref{2D_compare_densitydep}.

\subsection{Set-up of 2D radially-symmetric simulations}
\label{app:numerics_setup_radial}
To obtain the results displayed in Figure~\ref{2D_radial}, we solved numerically a 2D radially symmetric version of the PDE system~\eqref{nondim_PDE} with radial coordinate $r \in (0,10)$ and $t\in(0,10000]$, complemented with definition~\eqref{nondim_defineG} and subject to boundary conditions~\eqref{boundary}. Unless stated otherwise in the article, we set $L=10$. Numerical simulations were carried out using the {\sc{Matlab}} function {\bf{pdepe}}, with uniform discretisation of the spatial domain with step $\Delta_r=5\times10^{-3}$. 
As an initial condition, we chose the following perturbed version of the uniform steady state~\eqref{steadystate_app_switch}
\begin{equation}
{\small
n_0(0,r)=\overline{n}, \quad n_1(0,r)=1-\overline{n}, \quad s(0,r) = \overline{n}+0.01 \ e^{ -r^2},
\label{initialcondition_2Drs}
}
\end{equation}
where $\overline{n} = 0.5$. We note that, under the initial data defined via~\eqref{initialcondition_2Drs}, condition~\eqref{defsigmared} is satisfied. We set $D=1$ and $\chi=8$, and choose to investigate different values of $\mu$ as specified in the titles of the panels of Figure~\ref{2D_radial}.

\section{Numerical scheme employed in 2D simulations}
\label{app:2D_scheme}
The method for constructing 2D numerical solutions of the system~\eqref{nondim_PDE} complemented with definition~\eqref{nondim_defineG}, and subject to boundary conditions~\eqref{boundary} and condition~\eqref{defsigmared} is based on a modified version of the finite volume scheme employed in~\cite{bubba2020discrete}. The discretised dependent variables are
\begin{equation*}
	u^k_{i,j} := n_0(x_{i}, y_{j}, t_k), \quad v^k_{i,j} := n_1(x_{i}, y_{j}, t_k), \quad \text{and} \quad w^k_{i,j} := s(x_{i}, y_{j}, t_k).
\end{equation*}
First, we solve numerically the equation for the chemoattractant concentration, $w$, using the following scheme
\begin{equation*}
\begin{split}
\frac{w^{k+1}_{i,j}-w^{k}_{i,j}}{\tau} &= \frac{w^{k}_{i+1,j}-2w^{k}_{i,j} + w^{k}_{i-1,j}}{\left(\Delta _x\right)^2} +  \frac{w^{k}_{i,j+1}-2w^{k}_{i,j}+w^{k}_{i,j-1}}{\left(\Delta_y\right)^2} \\
						&\quad + u^k_{i,j} - w^k_{i,j}, \quad  i,j = 1, \dots, \mathcal{N},
\end{split}
\end{equation*}
and impose zero-flux boundary conditions by letting
\begin{equation*}
	\begin{split}
	w^{k+1}_{0,j} = w^{k+1}_{1,j}, \quad w^{k+1}_{\mathcal{N}+1,j} = w^{k+1}_{\mathcal{N},j}, \quad j = 1, \dots, \mathcal{N},\\
	w^{k+1}_{i,0} = w^{k+1}_{i,1}, \quad w^{k+1}_{i,\mathcal{N}+1} = w^{k+1}_{i,\mathcal{N}}, \quad i = 1, \dots, \mathcal{N}.
	\end{split}
\end{equation*}

Similarly, we solve numerically the equation for the density of cells in the non-chemotactic phenotypic state 0, $u$, using the following scheme
\begin{equation*}
\begin{split}
\frac{u^{k+1}_{i,j}-u^{k}_{i,j}}{\tau} &= D \frac{u^{k}_{i+1,j}-2u^{k}_{i,j} + u^{k}_{i-1,j}}{\left(\Delta _x\right)^2} + D \frac{u^{k}_{i,j+1}-2u^{k}_{i,j}+u^{k}_{i,j-1}}{\left(\Delta_y\right)^2} \\
						&\quad - \mu_{01}(u^k_{i,j},v^k_{i,j},w^k_{i,j})  u^k_{i,j} + \mu_{10}(u^k_{i,j},v^k_{i,j},w^k_{i,j})  v^{k}_{i,j}, \quad  i,j = 1, \dots, \mathcal{N},
\end{split}
\end{equation*}
and impose zero-flux boundary conditions by letting
\begin{equation*}
	\begin{split}
	u^{k+1}_{0,j} = u^{k+1}_{1,j}, \quad u^{k+1}_{\mathcal{N}+1,j} = u^{k+1}_{\mathcal{N},j}, \quad j = 1, \dots, \mathcal{N},\\
	u^{k+1}_{i,0} = u^{k+1}_{i,1}, \quad u^{k+1}_{i,\mathcal{N}+1} = u^{k+1}_{i,\mathcal{N}}, \quad i = 1, \dots, \mathcal{N}.
	\end{split}
\end{equation*} 

Then, we solve numerically the equation for the density of cells in the chemotactic phenotypic state 1, $v$, using the following scheme
\begin{equation*}
\begin{split}
	\frac{v^{k+1}_{i,j} - v^k_{i,j}}{\tau} =& \frac{F^{k}_{i+\frac{1}{2},j}- F^{k}_{i-\frac{1}{2},j}}{\Delta_x} + \frac{F^{k}_{i,j+\frac{1}{2}}-F^{k}_{i,j-\frac{1}{2}}}{\Delta_y}\\
	&+ \mu_{01}(u^k_{i,j},v^k_{i,j},w^k_{i,j})  u^k_{i,j} - \mu_{10}(u^k_{i,j},v^k_{i,j},w^k_{i,j})  v^{k}_{i,j}, \quad i, j = 1, \dots, \mathcal{N},
	\end{split}
\end{equation*}
where 
{\small
\begin{equation*}
\begin{split}
F^{k}_{i+\frac{1}{2},j} &= D \frac{v^k_{i+1,j}-v^k_{i,j}}{\Delta_x} - b^{k,+}_{i+\frac{1}{2},j} v^k_{i,j} +b^{k,-}_{i+\frac{1}{2},j} v^k_{i+1,j}, \quad i = 1,\dots,\mathcal{N}-1, \,\, j = 1, \dots, \mathcal{N},\\
F^{k}_{i,j+\frac{1}{2}} &= D \frac{v^k_{i,j+1}-v^k_{i,j}}{\Delta_y} - b^{k,+}_{i,j+\frac{1}{2}} v^k_{i,j}   +b^{k,-}_{i,j+\frac{1}{2}} v^k_{i,j+1}, \quad i = 1,\dots,\mathcal{N}, \,\, j = 1, \dots, \mathcal{N}-1,
\end{split}
\end{equation*}
}

with
{\small
\begin{equation*}
b^{k}_{i+\frac{1}{2},j} = \chi \frac{w^k_{i+1,j}-w^k_{i,j}}{\Delta_x}, \quad 	b^{k,+}_{i+\frac{1}{2},j} = \max \left(0, b^{k}_{i+\frac{1}{2},j} \right), \quad b^{k,-}_{i+\frac{1}{2},j} = \max \left(0, -b^{k}_{i+\frac{1}{2},j} \right),
\end{equation*}
}

and
{\small
\begin{equation*}
b^{k}_{i,j+\frac{1}{2}} = \chi \frac{w^n_{i,j+1}-w^n_{i,j}}{\Delta_y}, \quad b^{k,+}_{i,j+\frac{1}{2}} = \max \left(0, b^{k}_{i,j+\frac{1}{2}} \right), \quad b^{k,-}_{i,j+\frac{1}{2}} = \max \left(0, -b^{k}_{i,j+\frac{1}{2}} \right).
\end{equation*}
}

The discrete fluxes $F^{k}_{i-\frac{1}{2},j}$ for $i = 2,\dots,\mathcal{N}$, $j = 1, \dots, \mathcal{N}$ and $F^{k}_{i,j-\frac{1}{2}}$ for $i = 1,\dots,\mathcal{N}$, $j = 2, \dots, \mathcal{N}$ are defined in analogous ways, and we impose zero-flux boundary conditions by using the definitions
\begin{equation*}
\begin{split}
F^{k}_{1-\frac{1}{2},j} := 0, \quad F^{k}_{\mathcal{N}+\frac{1}{2},j} := 0, \quad j = 1, \dots, \mathcal{N},\\
F^{k}_{i,1-\frac{1}{2}} := 0,  \quad F^{k}_{i,\mathcal{N}+\frac{1}{2}} := 0, \quad i = 1, \dots, \mathcal{N}.
\end{split}
\end{equation*}
Notice that we employ a fully explicit scheme to avoid Newton sub-iterations that could be computationally expensive.  {{Computational times of simulations could be reduced by using more efficient numerical schemes, such as implicit-explicit schemes whereby diffusion terms are discretised implicitly but all other terms are discretised explicitly~\cite{hundsdorfer2003numerical}.}}

\section{Further figures for the cases where we observe non-stationary spatial patterns}

\begin{figure}[H]
\centering
\includegraphics[width=0.9\textwidth]{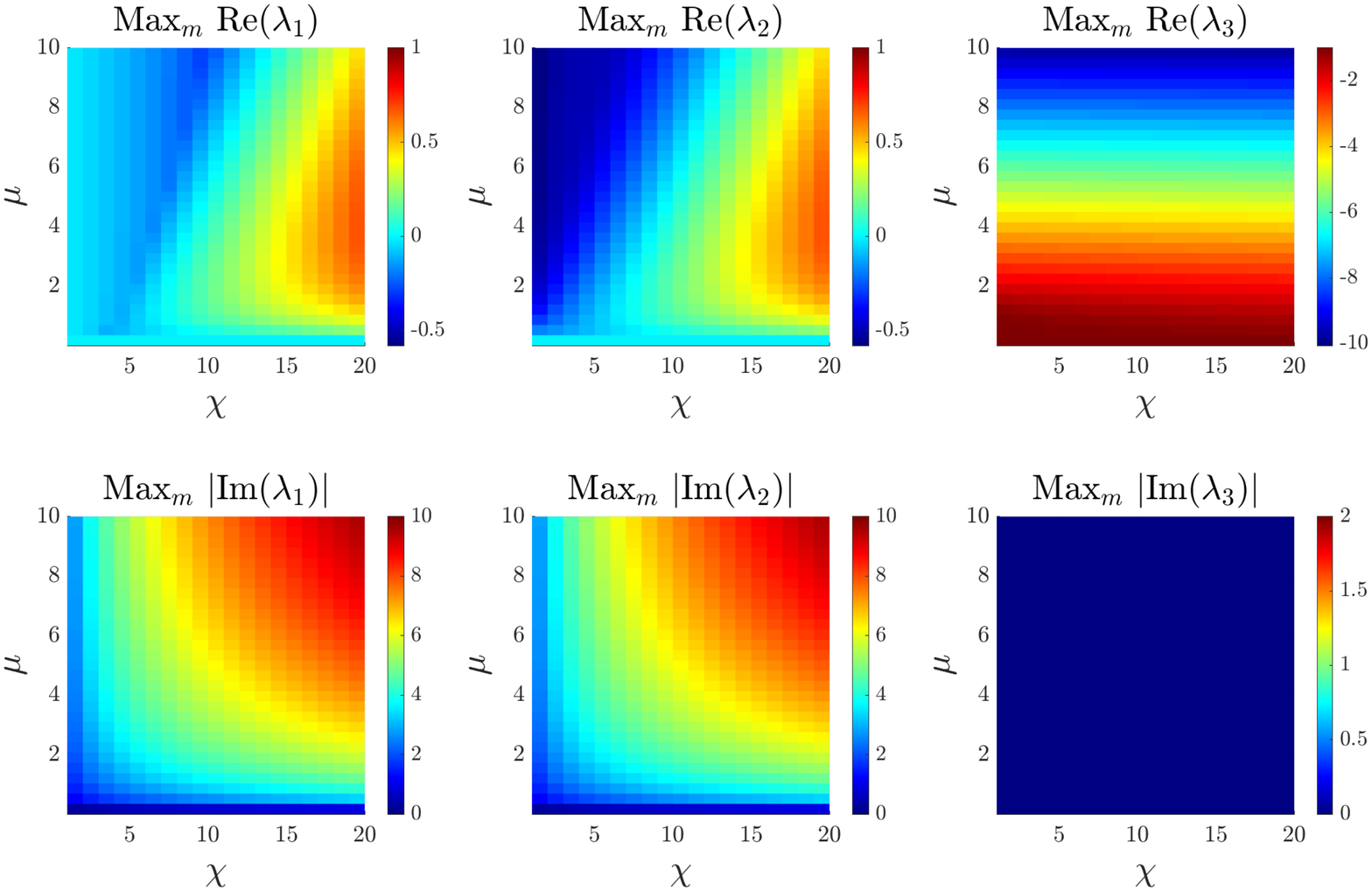}
\caption{{\bf{Investigating the eigenvalues for cases where we observe non-stationary spatial patterns.}} Here, the phenotypic switching functions are defined as Case B$_1$ given in Table~\ref{tab:1} and the parameter setting is $L=40$, $\overline{n}=0.5$, $D=1$, $q=30$. We compare the eigenvalues, i.e. the roots of~\eqref{characteristic_EQ}, for values of $\chi\in[1,20]$ and $\mu\in[0.01,10]$.  {\bf{Top row:}} We display the maximum value across all perturbation modes $m$ of the real part of each eigenvalue. {\bf{Bottom row:}} We display the maximum absolute value across all perturbation modes $m$ of the imaginary part of each eigenvalue.  }
\label{eigenvalues_oscillations}
\end{figure}

\bibliographystyle{plain}
\bibliography{bibliography}

\end{document}